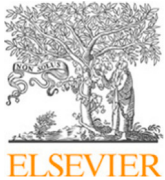
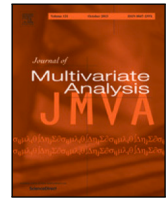

# Ledoit-Wolf linear shrinkage with unknown mean

Benoît Oriol [a,b], *, Alexandre Miot [b]

[a] *CEREMADE, Université Paris-Dauphine, PSL, place du Marechal de Lattre de Tassigny, 75116 Paris, France*
[b] *Société Générale Corporate and Investment Banking, 92800 Puteaux, France*

**ABSTRACT**

This work addresses large dimensional covariance matrix estimation with unknown mean. The empirical covariance estimator fails when dimension and number of samples are proportional and tend to infinity, settings known as Kolmogorov asymptotics. When the mean is known, Ledoit and Wolf (2004) proposed a linear shrinkage estimator and proved its convergence under those asymptotics. To the best of our knowledge, no formal proof has been proposed when the mean is unknown. To address this issue, we propose to extend the linear shrinkage and its convergence properties to translation-invariant estimators. We expose four estimators respecting those conditions, proving their properties. Finally, we show empirically that a new estimator we propose outperforms other standard estimators.

## 1. Introduction and related work

The covariance matrix plays a major role in numerous machine learning algorithms and statistics. Just to cite a few, the PCA [1] in machine learning, Markowitz portfolio management [2] in finance, or generalized method of moments estimators [3] in statistics. However, those algorithms are designed to use the true covariance matrix, which is often unaccessible. Even if the sample covariance matrix seems to be a simple and appealing choice, it severely fails in many applications: for instance, the use of the sample covariance matrix for Markowitz portfolio management does not beat a naive uniform distribution among the assets [4].

In the context of Kolmogorov asymptotics, where the ratio of the dimension $p_n$ and the number of samples $n$ tends to a finite positive constant $p_n/n \to c > 0$, this estimator fails to converge quadratically. Moreover, its eigenvalue spectrum is biased: high eigenvalues tend at being too high, and low ones, too low. The behavior of the eigenvalues is studied in random matrix theory: in the context of the Kolmogorov asymptotics, this topic is widely covered by V. L. Girko [5–7].

We focus on the shrinkage-type estimators which have suitable asymptotic properties, influenced by the work of Stein on Gaussian mean estimation in 1956 [8]. Due to their simplicity to implement and strong theoretical support, linear methods are widely used, and, for some, implemented in ScikitLearn [9]: Ledoit-Wolf linear shrinkage [10], which will be our main focus, its extension for Gaussian distributions using Rao–Blackwell theorem, named Oracle Approximating Shrinkage (OAS) estimator [11], linear shrinkage with factor models [12], linear shrinkage for elliptical distributions with unknown mean and known radius distribution [13], just to name a few. Non-linear methods propose shrinkage methods where the shrinkage intensity differs from an eigenvalue to an other. Among them, Stein's covariance estimator [14] works for Gaussian distributions, and several algorithms were developed by Ledoit and Wolf using eigenvalue spectrum analysis from random matrix theory [15–17]. Further theoretical analysis of those algorithms can be found in [18–20].

Usually, when estimating the covariance matrix, we do not know the mean of the distribution. Yet, the extension from known to unknown mean is rarely studied. To extend the empirical covariance with $T$ samples $S_T$, one uses the unbiased estimator $\tilde{S}_T$, after removing the empirical mean and dividing by $T-1$ instead of $T$. If it seems straightforward for $S_T$, it can be non trivial for more complex estimators. Ashurbekova explicitly worked in the case of known then unknown mean [13], and the resulting estimators of linear shrinkage for elliptical distributions with known radius distributions are notably different from the ones with a known mean. In the review of their work in 2020 [21], Ledoit and Wolf proved their results in the case where the mean is known, and they claim at the end "One then simply replaces $S_T$ with $\tilde{S}_T$ and $T$ with $T-1$ in all the previous descriptions and computations in practice" (Section 6: Computational Aspects and Code).

---

* Corresponding author at: CEREMADE, Université Paris-Dauphine, PSL, place du Marechal de Lattre de Tassigny, 75116 Paris, France.
*E-mail address:* benoit.oriol@dauphine.eu (B. Oriol).

Theoretically, Ledoit and Wolf [17] used an interesting result in Remark 2.1, taken from the beginning of section 3 of Silverstein and Bai [22], to show that the perturbation caused by the estimation of the mean does not change the asymptotic spectrum of the sample covariance. Under the hypotheses of the paper [17], the asymptotic non-linear shrinkage oracle formulas depend only on this spectrum, and this lemma leads to show that those oracle formulas are the same if we exactly know or if we estimate the mean of the distribution.

However, while it works to justify the similar behavior of both oracle non-linear shrinkage formulas, when mean is known or not, under suitable hypotheses, it says nothing about the estimators of those oracles, which are used in practice. Regarding the linear shrinkage intensities, the estimators are not described as a function of the spectrum, their convergence properties cannot be derived directly from the result from Silverstein and Bai, and a deeper work has to be done to extend the theoretical properties.

Moreover, focusing on Ledoit-Wolf linear shrinkage algorithm, one can note that the implementation used in ScikitLearn [9] does not follow the recommendations of Ledoit and Wolf regarding the case of uncentered data. They did not change $T$ to $T-1$ and used $(T-1)\tilde{S}_T/T$ instead of $\tilde{S}_T$. Unexpectedly, experiments show notably worse results using Ledoit and Wolf recommendations rather than the ScikitLearn implementation. This remark underlines that the problem is more counter-intuitive than expected, and a closer look at the dependence between the covariance and the mean estimation is required.

The question this work answers is: how to estimate the oracle linear shrinkage intensities when the mean is unknown, and, do we retrieve the same convergence properties we had in the more simple case where the mean is known ? Two candidates of estimators are proposed by Ledoit Wolf in 2019 and by the implementation in Scikit Learn, and two new ones emerge naturally from the theoretical part.

We address the lack of theoretical results when the mean is unknown, extending the theoretical properties of linear shrinkage to four translation-invariant estimators, including the one currently implemented one in ScikitLearn 1.2.2. The experimental part aims at giving an understanding of which estimator to use in different finite samples concentrations and distributions, compare their empirical convergence speed, and confirm that the experimental robustness highlighted by Ledoit and Wolf [10] is conserved.

## 2. Notation, definitions and hypotheses

Let us introduce the following notation. To remain consistent with the original work of Ledoit and Wolf [10], we will use the same notation, conventions and assumptions.

In order to study theoretical properties under general asymptotics where the dimension and the number of samples are not constant, we consider a sequence of observations matrices $(X_n)$, indexed by $n$, the number of samples.

**Notation 1.** *In the following we consider a sequence of observation matrices $(X_n)_{n \in \mathbb{N}^*}$ with $X_n \in \mathbb{R}^{p_n \times n}$ of $n$ iid observations on a system of $p_n$ dimensions. Decomposing the covariance matrix, we denote $\Sigma_n = \Gamma_n \Lambda_n \Gamma_n^\top$, where $\Lambda_n$ is a diagonal matrix and $\Gamma_n$ a rotation matrix. The diagonal elements of $\Lambda_n$ are the eigenvalues $\lambda_1^n, \ldots, \lambda_{p_n}^n$, and the columns of $\Gamma_n$ are the eigenvectors $\gamma_1^n, \ldots, \gamma_{p_n}^n$. $Y_n = \Gamma_n^\top X_n$ is a $p_n \times n$ matrix of $n$ iid observations of $p_n$ uncorrelated random variables $(y_1^n, \ldots, y_n^n)$.*

As the dimension of the spaces we consider along with the sequence of observations is not constant, the choice of a norm for the convergence has an impact. Here, all norms are not equivalent. Even if at each step we are in a space of matrices of finite dimension $(p_n, p_n)$, this dimension evolves with $n$, and can diverge. Thus we chose the Frobenius norm, normalized by $\sqrt{p_n}$. It is not standard, but it aims at giving the norm 1 to the identity $I_{p_n}$ regardless of the dimension $(p_n, p_n)$.

**Notation 2.** *Let $A_n$ and $B_n$ two $p_n \times p_n$ matrices. We consider the Frobenius norm: $\|A_n\| = \sqrt{\mathrm{tr}(A_n A_n^\top)/p_n}$, and the associated inner product: $\langle A_n, B_n \rangle_n = \mathrm{tr}(A_n B_n^\top)/p_n$. Dividing by the dimension is not standard, it is done to fix the norm of the identity as 1 regardless of the dimension.*

We naturally define what we mean by quadratic convergence of a sequence of random variables with non-constant dimensions.

**Notation 3.** *Let $(E_n)_n$ a sequence of euclidean spaces with associated norm $\|\cdot\|_n$. The quadratic convergence of a random variable $Z_n \in E_n$, i.e., $\mathbb{E}[\|Z_n\|_n^2] \to 0$, is denoted as $Z_n \xrightarrow[q.m]{} 0$.*

We describe now several assumptions, the same used in the linear shrinkage of Ledoit and Wolf [10], that will be used in the following. The first one defines the scope of the asymptotics. This is a bit less restrictive than the usual definition of Kolmogorov asymptotics, that imposes that $p_n/n \to c \in \mathbb{R}_+$. We only requires the quotient to be bounded here, the convergence is not necessary.

**Assumption 1.** *There exists a constant $K_1$ independent of $n$ such that $p_n/n \leq K_1$.*

The second assumption puts a condition on the eighth moment of the distributions, that is to be bounded in average along the dimensions. This assumption is quite strong and we need it for the proofs. However, experimentally, it seems that fourth moments are sufficient for the convergence of the estimators. Ledoit and Wolf [10] denoted it in the known mean case too.

**Assumption 2.** There exists a constant $K_2$ independent of $n$ such that $\frac{1}{p_n}\sum_{i=1}^{p_n} \mathbb{E}\left[(\tilde{y}_{i1}^n)^8\right] \leq K_2$ where for all $i \in [\![1, p_n]\!]$, $\tilde{y}_{i1}^n = y_{i1}^n - \mathbb{E}[y_{i1}^n]$.

This last assumption is technical, and Ledoit and Wolf [10] in section 3.1 ensured the usual distributions respect the condition, such as random variables normal or elliptical distributions, but it is much weaker than that.

**Assumption 3.**
$$\lim_{n\to\infty} \frac{p_n^2}{n} \times \frac{\sum_{(i,j,k,l)\in Q_n}(\mathrm{Cov}[\tilde{y}_{i1}^n \tilde{y}_{j1}^n, \tilde{y}_{k1}^n \tilde{y}_{l1}^n])^2}{|Q_n|} = 0,$$
where $Q_n$ denotes the set of all the quadruples that are made of four distinct integers between 1 and $p_n$, and for all $i \in [\![1, p_n]\!]$, $\tilde{y}_{i1}^n = y_{i1}^n - \mathbb{E}[y_{i1}^n]$.

We need some definitions to properly define the problem and the asymptotics. The main object of interest is the empirical covariance, or the sample covariance.

**Definition 1** (*Empirical Covariance*). For an observation matrix $X_n$ of size $p_n \times n$, we define the empirical covariance as:
$$S_n = \tilde{X}_n \tilde{X}_n^\top / (n-1),$$
with $(\tilde{X}_n)_{ik} = (X_n)_{ik} - \frac{1}{n}\sum_{k'=1}^n (X_n)_{ik'}$.

We introduce four scalars that will help us to characterize the oracle linear shrinkage.

**Definition 2** (*Scalars $(\mu_n, \alpha_n^2, \beta_n^2, \delta_n^2)$*). We define four scalars:
$$\mu_n = \langle \Sigma_n, I_{p_n} \rangle_n, \alpha_n^2 = \|\Sigma_n - \mu_n I_{p_n}\|_n^2, \beta_n^2 = \mathbb{E}[\|S_n - \Sigma_n\|_n^2], \delta_n^2 = \mathbb{E}[\|S_n - \mu_n I_{p_n}\|_n^2].$$
(Lemma 2.1 in [10]) proves that $\alpha_n^2 + \beta_n^2 = \delta_n^2$.

The oracle linear shrinkage estimator is given by the following minimization problem. The following corollary is the central point of the linear shrinkage methods.

**Corollary 1** (*Corollary of Theorem 2.1 from [10]*). *Consider the optimization problem:*
$$\begin{array}{c}\underset{\rho_1,\rho_2}{\text{minimize}} \quad \mathbb{E}[\|\Sigma_n^* - \Sigma_n\|_n^2], \\ \text{s.t.} \quad \Sigma_n^* = \rho_1 I_{p_n} + \rho_2 S_n\end{array}$$
where the coefficients $\rho_1$ and $\rho_2$ are not random. Its solution $\Sigma_n^*$ verifies:
$$\Sigma_n^* = \frac{\beta_n^2}{\delta_n^2}\mu_n I_{p_n} + \frac{\alpha_n^2}{\delta_n^2} S_n, \quad \mathbb{E}[\|\Sigma_n^* - \Sigma_n\|_n^2] = \frac{\alpha_n^2 \beta_n^2}{\delta_n^2}.$$

**Remark 1.** Corollary 1 remains true for any unbiased estimator $\hat{\Sigma}_n$ instead of $S_n$.

$(\mu_n, \alpha_n^2, \beta_n^2, \delta_n^2)$ depends on the true covariance $\Sigma_n$, and thus cannot be used directly in the estimation of $\Sigma_n^*$. The central issue of this work is to find estimators $(m_n, a_n^2, b_n^2, d_n^2)$ of $(\mu_n, \alpha_n^2, \beta_n^2, \delta_n^2)$ in order to compute an estimation $S_n^*$ of $\Sigma_n^*$. As the mean is unknown, those estimators differ from Ledoit and Wolf work [10], particularly when $p_n$ is higher than $n$.

## 3. Theoretical results

All the following results extend the work of Ledoit and Wolf [10] in the case where the empirical mean is used as estimator of the mean.

All proofs are shown in Appendix.

**Remark 2.** In the following, as all the estimators are invariant by change of mean, resulting from the definition of $\tilde{X} = X - \sum_k X_{\cdot,k}$, we can assume $\mathbb{E}[X] = 0$ for the simplicity of notation.

We present a sequence of lemmata, that naturally define estimators with suitable asymptotic properties for the scalars $(\mu_n, \alpha_n^2, \beta_n^2, \delta_n^2)$.

**Lemma 1.** *Under Assumptions 1 and 2, $\mu_n, \alpha_n^2, \beta_n^2, \delta_n^2$ remain bounded as $n \to \infty$.*

**Theorem 1.** *Under Assumptions 1 and 2, define $\theta_n^2 = \mathbb{V}\left[\frac{1}{p_n}\sum_{i=1}^{p_n}(y_{i1}^n)^2\right]$. $\theta_n^2$ is bounded as $n \to \infty$, and we have:*
$$\lim_{n\to\infty} \mathbb{E}[\|S_n - \Sigma_n\|_n^2] - \frac{p_n}{n}(\mu_n^2 + \theta_n^2) = 0.$$

In particular, taking $p_n = p$ constant, we see that $S_n \xrightarrow[q.m]{} \Sigma_n$. But, when $p_n$ is of the same order of magnitude than $n$, the sample covariance generally fails to converge as the error is at least of the same order of magnitude as $\mu_n^2 = \langle \Sigma_n, I_{p_n}\rangle_n^2$.

The next lemma defines a natural estimator of $\mu_n$, and its convergence properties are conserved when the mean is estimated.

**Lemma 2** (*Estimator of $\mu_n$*). *Define $m_n = \langle S_n, I_{p_n}\rangle_n$. Then, under Assumptions 1 and 2, $\mathbb{E}[m_n] = \mu_n$ for all $n$, and $m_n - \mu_n$ converges to zero in quartic mean (fourth moment) as $n$ goes to infinity.*

**Corollary 2.** *Under Assumptions 1 and 2, $m_n^2 - \mathbb{E}[m_n^2]$ converges to zero in quadratic mean as $n$ goes to infinity.*

For $\delta_n^2$, the natural estimator $d_n^2$ defined below, and its convergence properties are conserved when the mean is estimated.

**Lemma 3** (*Estimator of $\delta_n^2$*). *Define $d_n^2 = \|S_n - m_n I_{p_n}\|_n^2$. Then, under Assumptions 1, 2 and 3, $d_n^2 - \delta_n^2 \xrightarrow[q.m]{} 0$. It follows that $d_n^2 - \mathbb{E}[d_n^2] \xrightarrow[q.m]{} 0$.*

The main difference that we have when estimating the mean is around the estimation of $\beta_n^2$. The following lemmata aim at defining an unbiased estimator of $\beta_n^2$ that quadratically converges. We work around $\bar{b}_n^2$, inspired from the estimator of $\beta_n^2$ in the case where the mean is known, and write residual terms in the expectation as a combination of $m_n$ and $d_n^2$.

**Lemma 4.** *Define:*
$$\bar{b}_n^2 = \frac{1}{n^2}\sum_{k=1}^n \left\|\frac{n}{n-1}\tilde{x}_{\cdot k}^n(\tilde{x}_{\cdot k}^n)^t - S_n\right\|_n^2,$$
*where $\tilde{x}_{\cdot k}^n = x_{\cdot k}^n - \frac{1}{n}\sum_{k'=1}^n x_{\cdot k'}^n$, and $(x_{\cdot 1}^n, \ldots, x_{\cdot k}^n, \ldots, x_{\cdot n}^n)$ are the independent samples forming $X$. Then, under Assumption 1,*
$$\mathbb{E}[\bar{b}_n^2] = c_0 \beta_n^2 + c_1 \delta_n^2 + c_2 \mu_n^2,$$
*with $\gamma_n = \frac{n(n-1)}{n^2 - 3n + 3}$, $\lambda_n = \frac{n^2(n-2)}{(n-1)(n^2 - 3n + 3)}$, $c_0 = \frac{1}{\gamma_n} - \frac{1}{n} - \frac{\lambda_n}{\gamma_n n^2}$, $c_1 = \frac{\lambda_n}{\gamma_n n^2}$, $c_2 = (p+1)c_1$.*

**Lemma 5.** *Under Assumption 1, we have: $\mathbb{E}[\bar{b}_n^2] = c_0 \beta_n^2 + c_1 \mathbb{E}[d_n^2] + c_2 \mathbb{E}[m_n^2] + (c_1 - c_2)\mathbb{V}[m_n]$.*

**Lemma 6.** *Under Assumption 1, we have: $\mathbb{V}[\bar{b}_n^2] \xrightarrow[n\to\infty]{} 0$.*

We need to compute $\mathbb{V}[m_n]$ which is unknown in the development of $\mathbb{E}[\bar{b}_n^2]$ in Lemma 5. We want to express $\mathbb{V}[m_n]$ as a combination of expectations of quantities we can compute with $X_n$, as $m_n$ and $d_n^2$, so that we can find an unbiased estimator of $\beta_n^2$ using Lemma 5.

**Lemma 7.** *Under Assumption 1, we have: $\mathbb{V}[m_n] = q_0 \beta_n^2 + q_1 \delta_n^2 - q_2 \mu_n^2$, with $q_0 = \frac{n-2}{p(n-1)}$, $q_1 = \frac{1}{p(n-1)}$, $q_2 = \frac{p-1}{p(n-1)}$.*

In this first step, we write $\mathbb{V}[m_n]$ as a combination of $\beta_n^2, \delta_n^2, \mu_n$. We then use the expressions of $\mathbb{E}[d_n^2]$ and $\mathbb{E}[m_n^2]$ to finish the work about $\mathbb{V}[m_n]$.

**Lemma 8.** *Under Assumption 1, we have: $\mathbb{V}[m_n] = \frac{1}{1-q_1-q_2}\left(q_0 \beta_n^2 + q_1 \mathbb{E}[d_n^2] - q_2 \mathbb{E}[m_n^2]\right)$.*

Backing up, using $\bar{b}_n^2$, $d_n^2$ and $m_n$, we define $b_n^2$ an unbiased estimator of $\beta_n^2$ which converges in quadratic mean to $\beta_n^2$.

**Lemma 9.** *Define:*
$$b_n^2 = \frac{1}{c_0^f}\left(\bar{b}_n^2 - c_1^f d_n^2 - c_2^f m_n^2\right),$$
*with $c_0^f = c_0 + (c_1 - c_2)\frac{q_0}{1-q_1-q_2}$, $c_1^f = c_1 + (c_1 - c_2)\frac{q_1}{1-q_1-q_2}$, $c_2^f = c_2 - (c_1 - c_2)\frac{q_2}{1-q_1-q_2}$.*
*Then, under Assumptions 1, 2 and 3, $b_n^2$ is an unbiased estimator of $\beta_n^2$, i.e. $\mathbb{E}[b_n^2] = \beta_n^2$, and $b_n^2 - \beta_n^2 \xrightarrow[q.m]{} 0$.*

For notation consistency with the estimators in Ledoit-Wolf linear shrinkage [10], we keep the notation $b_n^2$ even if its value can be negative.

As in Ledoit and Wolf work [10], the final estimator $b_{n,u}^2$ of $\beta_n^2$ goes through a threshold so that the estimated shrinkage intensity $b_{n,u}^2/d_n^2 \in [0,1]$. It does not change the theoretical properties but experimentally increased significantly the results.

**Lemma 10** (*Estimator of $\beta_n^2$*). *Define: $b_{n,u}^2 = \min((b_n^2)_+, d_n^2)$ and $a_{n,u}^2 = d_n^2 - b_{n,u}^2$. Under Assumptions 1, 2 and 3, $b_{n,u}^2 - \beta_n^2 \xrightarrow[q.m]{} 0$ and $a_{n,u}^2 - \alpha_n^2 \xrightarrow[q.m]{} 0$.*

We can now define our translation-invariant linear shrinkage estimator $S_n^*$ and prove its asymptotic properties.

**Definition 3** (*Final Estimator LW_u*). Let us define our estimator:

$$S_n^* = \frac{b_{n,u}^2}{d_{n,u}^2} m_{n,u} I_{p_n} + \frac{a_{n,u}^2}{d_{n,u}^2} S_n, \text{ with } m_{n,u} = m_n, d_{n,u}^2 = d_n^2.$$

We recover the same quadratic convergence to $\Sigma_n^*$ with our translation-invariant estimator as we had originally when the mean is known.

**Theorem 2.** *Under Assumptions 1, 2 and 3, $\mathbb{E}[\|S_n^* - \Sigma_n^*\|^2] \to 0$. As a consequence, $S_n^*$ has the same asymptotic expected loss as $\Sigma_n^*$, i.e. $\mathbb{E}[\|S_n^* - \Sigma_n\|_n^2] - \mathbb{E}[\|\Sigma_n^* - \Sigma_n\|_n^2] \to 0$.*

The following lemma gives an asymptotic estimation of the optimal error $\mathbb{E}\left[\|\Sigma_n^* - \Sigma_n\|_n^2\right] = \frac{\alpha_n^2 \beta_n^2}{\delta_n^2}$.

**Lemma 11.**

$$\mathbb{E}\left[\left|\frac{a_{n,u}^2 b_{n,u}^2}{d_n^2} - \frac{\alpha_n^2 \beta_n^2}{\delta_n^2}\right|\right] \to 0.$$

The last results easily make possible to extend the Theorems 3.3 and 3.4 of (Ledoit and Wolf, 2004) [10] in our situation where the mean is unknown. Previously, we showed that our estimator's loss converge to the optimal one in the class of linear combinations of $S_n$ and $I_{p_n}$ with non random coefficients, the optimal estimator of this class being $\Sigma_n^*$. In the following, we show that our estimator is still asymptotically optimal with respect to a bigger class, where the coefficients can be random. Formally, we are looking for the following optimal loss (this time, there is no expectation in the minimization). Let $\Sigma_n^{**}$ be the linear combination of $S_n$ and $I_{p_n}$ solving:

$$\underset{\rho_1, \rho_2}{\text{minimize}} \quad \|\Sigma_n^{**} - \Sigma_n\|_n^2.$$
$$\text{s.t. } \Sigma_n^{**} = \rho_1 I_{p_n} + \rho_2 S_n$$

By construction, $\Sigma_n^{**}$ has a lower loss than $S_n^*$, but we show that the difference converges to 0.

**Theorem 3.** *$S_n^*$ converges to $\Sigma_n^{**}$ in quadratic mean, i.e. $\|S_n^* - \Sigma_n^{**}\| \underset{q.m}{\longrightarrow} 0$. As a consequence, $S_n^*$ has the same asymptotic expected loss as $\Sigma_n^{**}$, more precisely we have:*

$$\mathbb{E}\left[\left|\|S_n^* - \Sigma_n\|_n^2 - \|\Sigma_n^{**} - \Sigma_n\|_n^2\right|\right] \to 0.$$

**Theorem 4.** *For any sequence of linear combinations $\hat{\Sigma}_n$ of $I_n$ and $S_n$, the estimator $S_n^*$ verifies:*

$$\lim_{N \to \infty} \inf_{n \geq N} \left(\mathbb{E}\left[\|\hat{\Sigma}_n - \Sigma_n\|_n^2\right] - \mathbb{E}\left[\|S_n^* - \Sigma_n\|_n^2\right]\right) \geq 0.$$

*In addition, every $\hat{\Sigma}_n$ that performs as well as $S_n^*$ is identical to $S_n^*$ in the limit:*

$$\lim_{N \to \infty} \left(\mathbb{E}\left[\|\hat{\Sigma}_n - \Sigma_n\|_n^2\right] - \mathbb{E}\left[\|S_n^* - \Sigma_n\|_n^2\right]\right) = 0 \iff \|\hat{\Sigma}_n - S_n^*\|_n \underset{q.m}{\longrightarrow} 0.$$

We introduce three other translation-invariant estimators to compare with, which are implemented, recommended, or natural to define. We prove that their asymptotic behavior is similar, and, through different experiments, show the differences in performance.

The first one is the estimator recommended by Ledoit and Wolf in their Section 6: Computational Aspects and Code [21].

**Definition 4** (*Ledoit-Wolf Recommended Estimators*). The estimators recommended by Ledoit and Wolf [21], indexed by the letter "$r$", are:

$$m_{n,r} = m_n, d_{n,r}^2 = d_n^2, b_{n,r}^2 = \min\left(\left(\frac{1}{(n-1)^2}\sum_{k=1}^n \left\|\tilde{x}_{\cdot k}^n(\tilde{x}_{\cdot k}^n)^t - S_n\right\|_n^2\right)_+, d_{n,r}^2\right), a_{n,r}^2 = d_{n,r}^2 - b_{n,r}^2.$$

The following theorem ensures that all the convergence results in Ledoit and Wolf [10] remain true with this translation-invariant estimator.

**Theorem 5** (*Ledoit-Wolf Recommended Estimators*). *Under Assumptions 1, 2 and 3, $m_{n,r} - \mu_n$ converges to 0 in quartic mean, and that $d_{n,r}^2 - \delta_n^2$, $b_{n,r}^2 - \beta_n^2$ and $a_{n,r}^2 - \alpha_n^2$ converge in quadratic mean to 0 as $n$ goes to infinity.*

*Moreover, the conclusions of Theorem 2 remain true with the estimated matrix $S_{n,r}^* = \frac{b_{n,r}^2}{d_{n,r}^2} m_{n,r} I_{p_n} + \frac{a_{n,r}^2}{d_{n,r}^2} S$, i.e. $\mathbb{E}[\|S_{n,r}^* - \Sigma_n^*\|^2] \to 0$ and*

$$\mathbb{E}\left[\left|\|S_{n,r}^* - \Sigma_n\|_n^2 - \|\Sigma_n^{**} - \Sigma_n\|_n^2\right|\right] \to 0.$$

From the proof, $\bar{b}_{n,r}^2 = \bar{b}_n^2 + \frac{1}{n(n-1)^2}\|S_n\|_n^2$, it is then natural to define the following estimator.

**Definition 5** ("*Natural*" *Estimators*). The estimators that naturally emerge, indexed by the letter "m", are:
$$m_{n,m} = m_n, d_{n,m}^2 = d_n^2, b_{n,m}^2 = \min\left(\left(\bar{b}_n^2\right)_+, d_{n,m}^2\right), a_{n,m}^2 = d_{n,m}^2 - b_{n,m}^2.$$

The following theorem ensures too that all the convergence results in Ledoit and Wolf [10] remain true with this natural translation-invariant estimator.

**Theorem 6** ("*Natural*" *Estimators*). *Under Assumptions 1, 2 and 3, $m_{n,m} - \mu_n$ converges to 0 in quartic mean, and that $d_{n,m}^2 - \delta_n^2$, $b_{n,m}^2 - \beta_n^2$ and $a_{n,m}^2 - \alpha_n^2$ converge in quadratic mean to 0 as $n$ goes to infinity.*

*Moreover, the conclusions of Theorem 2 remain true with the estimated matrix $S_{n,m}^* = \frac{b_{n,m}^2}{d_{n,m}^2} m_{n,m} I_{p_n} + \frac{a_{n,m}^2}{d_{n,m}^2} S$, i.e. $\mathbb{E}[\|S_{n,m}^* - \Sigma_n^*\|^2] \to 0$ and $\mathbb{E}\left[\left|\|S_{n,m}^* - \Sigma_n\|_n^2 - \|\Sigma_n^{**} - \Sigma_n\|_n^2\right|\right] \to 0$.*

The final estimator is the one implemented in ScikitLearn 1.2.2 when the mean is unknown. One goal of this theoretical part is to prove whether or not this estimator widely used in practice has indeed the theoretical properties we expected it to have from the known mean situation.

**Definition 6** (*ScikitLearn 1.2.2 Estimators*). The estimators implemented in ScikitLearn 1.2.2, indexed by the letter "s", are:
$$m_{n,s} = \frac{n-1}{n} m_n, d_{n,s}^2 = d_n^2, b_{n,s}^2 = \min\left(\left(\bar{b}_n^2\right)_+, d_{n,s}^2\right), a_{n,s}^2 = \frac{n-1}{n}\left(d_{n,s}^2 - b_{n,s}^2\right),$$

$$S_{n,s}^* = \frac{n-1}{n} S_{n,m}^*.$$

This final theorem ensures that all the convergence results in Ledoit and Wolf [10] remain true with the implementation of ScikitLearn 1.2.2 when the mean is estimated.

**Theorem 7** (*Scikit-Learn 1.2.2 Estimators*). *Under Assumptions 1, 2 and 3, $m_{n,s} - \mu_n$ converges to 0 in quartic mean, and that $d_{n,s}^2 - \delta_n^2$, $b_{n,s}^2 - \beta_n^2$ and $a_{n,s}^2 - \alpha_n^2$ converge in quadratic mean to 0 as $n$ goes to infinity.*

*Moreover, the conclusions of Theorem 2 remain true with the estimated matrix $S_{n,s}^* = \frac{b_{n,s}^2}{d_{n,s}^2} m_{n,s} I_{p_n} + \frac{a_{n,s}^2}{d_{n,s}^2} S$, i.e. $\mathbb{E}[\|S_{n,s}^* - \Sigma_n^*\|^2] \to 0$ and $\mathbb{E}\left[\left|\|S_{n,s}^* - \Sigma_n\|_n^2 - \|\Sigma_n^{**} - \Sigma_n\|_n^2\right|\right] \to 0$.*

## 4. Experimental results

As we proved it in the theoretical results, the four estimators we study have the same asymptotic properties, the same we expected when the mean is known. However, the speed of convergence, and thus performance at finite sample are *a priori* different between the estimators. Here, we show that most of the difference is at high $p_n/n$, when the proposed estimator $S_n^*$ converges significantly faster than the other candidates. In low concentration, the behaviors are very similar, as expected. Indeed, when we have a lot of samples compared to the dimension, we expect that the loss of information for estimating the mean is little in comparison, and thus any correction based on the variance of the estimation of mean may not change much the whole estimator.

The experimental estimations are compared to the theoretical value of $\Sigma^*$ and $\Sigma^{**}$ in the Ledoit-Wolf setting, the implementation in ScikitLearn 1.2.2, the implementation recommended by Ledoit and Wolf [21], and to the other algorithms implemented in ScikitLearn 1.2.2, for multivariate Gaussian and Student-t distributions.

We first derive the exact values of $\Sigma^*$ for those two distributions.

### 4.1. Oracle estimators

#### 4.1.1. Gaussian distribution
**Lemma 12.** *Let $(X_n)_{\cdot,k} \sim \mathcal{N}(0, \Sigma_n)$, $k \in [\![1,n]\!]$, $n$ iid samples. Then, the analytical oracle estimators are: $\mu_n = \langle \Sigma_n, I \rangle_n, \alpha_n^2 = \|\Sigma_n\|_n^2 - \mu_n^2, \beta_n^2 = \frac{p+1}{n-1} \mu_n^2 + \frac{1}{n-1} \alpha_n^2, \delta_n^2 = \alpha_n^2 + \beta_n^2$.*

#### 4.1.2. t-distribution oracle estimator
**Lemma 13.** *Let $(X_n)_{\cdot,k} \sim t_\nu(0, \tilde{\Sigma}_n)$, $k \in [\![1,n]\!]$, $n$ iid samples with scale matrix $\tilde{\Sigma}_n = \frac{\nu-2}{\nu} \Sigma_n$ and covariance $\mathbb{V}[X_n] = \Sigma_n$. The density of the multivariate t-distribution is:*
$$f_n(x) = \frac{\Gamma[(\nu+p)/2]}{\Gamma(\nu+p)\nu^{p/2}\pi^{p/2}|\tilde{\Sigma}_n|^{1/2}} \left[1 + \frac{1}{\nu} x^\top \tilde{\Sigma}_n^{-1} x\right]^{-(\nu+p)/2}.$$

*Then, the analytical oracle estimators are:*
$\mu_n = \langle \Sigma_n, I \rangle, \alpha_n^2 = \|\Sigma_n\|_n^2 - \mu_n^2, \beta_n^2 = \frac{1}{n}\left(\frac{\nu}{\nu-4} + \frac{1}{n-1}\right)\left(\alpha_n^2 + (p+1)\mu_n^2\right) - \frac{2p}{n(\nu-4)} \mu_n^2, \delta_n^2 = \alpha_n^2 + \beta_n^2$.

### 4.2. Experimental setup

We considered 2 settings:

- a Monte-Carlo computation of the loss on a 2d-grid of the parameters $(p_n, n) \in [\![5, 100]\!]^2$, with a step size of 2, to visualize the effect of changing the ratio $p_n/n$, and see the domains where our algorithm is most suited;
- a Monte-Carlo computation of the loss with a fixed ratio $p_n/n = c \in \{0.25, 0.5, 1, 2, 4\}$, to compare the rate of convergence of each algorithm.

In both cases, the Monte-Carlo is computed with $n_{MC} = 10000$ iterations.

Three types of distributions are explored: the multivariate Gaussian, the t-distribution with $\nu = 10$, $\nu = 8.5$ and $\nu = 4.5$, and a mix of two independent t-distributions with $\nu = 15$ on the first half coordinates, and $\nu = 8.5$ on the last half. Note that we have to ensure $\nu > 8$ to respect Assumption 2. While $\nu = 4.5$ is outside of the scope of the theoretical part, we want to show the experimental robustness of linear shrinkage to heavy tails, as it was noticed originally in Ledoit and Wolf [10].

Two different way of choosing $\Sigma_n$ are explored: fixing $\Sigma_n = I_{p_n}$ - particular case where the oracle Ledoit-Wolf loss is null -, and drawing at each iteration a covariance matrix $\Sigma_n$ from a Wishart distribution with $p_n$ degrees of freedom, and normalizing it by $\sqrt{\|\Sigma_n \Sigma_n^\top\|_n}$ - to respect Assumption 2. Note that when drawn from a Wishart, $\sqrt{\|\Sigma_n \Sigma_n^\top\|_n} > 0$ almost surely.

#### 4.2.1. Assumptions check

For the second study at fixed $p_n/n$ in order to compare the rate of convergence, we check that we are under the three assumptions that guarantee the theoretical results on convergence proved in Section 3.

As we fixed $p_n/n$, Assumption 1 is trivially respected.

*Assumption 2 - Gaussian distribution.* Let $(X_n)_{\cdot,k} \sim \mathcal{N}(0, \Sigma_n)$, $k \in [\![1, n]\!]$, $n$ iid samples. As previously, we denote $\Sigma_n = \Gamma_n \Lambda_n \Gamma_n^\top$, where $\Lambda_n$ is a diagonal matrix and $\Gamma_n$ a rotation matrix, and $Y_n = \Gamma_n^\top X_n$ is a $p_n \times n$ matrix of $n$ iid observations of $p_n$ uncorrelated random variables.

Using the fact that for $z \sim \mathcal{N}(0, \lambda)$, we have $\mathbb{E}[z^8] = 105 \lambda^4$, we deduce:

$$\frac{1}{p_n} \sum_{i=1}^{p_n} \mathbb{E}[y_{i1}^8] = \frac{1}{p_n} \sum_{i=1}^{p_n} 105 \lambda_{ii}^4 = 105 \|\Lambda_n \Lambda_n^\top\|_n^2 = 105 \|\Sigma_n \Sigma_n^\top\|_n^2. \tag{1}$$

In the case where we fix $\Sigma_n = I_{p_n}$, we obviously have $\|\Sigma_n \Sigma_n^\top\|_n^2 = 1$, so Assumption 2 is respected. In the case where we draw $\Sigma_n$ from a Wishart distribution with $p_n$ degrees of freedom, and normalize it by $\sqrt{\|\Sigma_n \Sigma_n^\top\|_n}$, we have by construction $\|\Sigma_n \Sigma_n^\top\|_n^2 = 1$, so Assumption 2 is respected here too.

*Assumption 2 - t-distribution.* Let $(X_n)_{\cdot,k} \sim t_\nu(0, \tilde{\Sigma}_n)$, $k \in [\![1, n]\!]$, $n$ iid samples with $\nu > 8$, scale matrix $\tilde{\Sigma}_n = \frac{\nu-2}{\nu} \Sigma_n$ and covariance $\mathbb{V}[X_n] = \Sigma_n$.

From a characterization of multivariate t-distributions, for each $k \in [\![1, n]\!]$, there exist 2 independent random variables $U_k$ and $Z_{\cdot,k}$ such that:

$$U_k \sim \chi_\nu^2, Z_{\cdot,k} \sim \mathcal{N}\left(0, \frac{\nu-2}{\nu}\Lambda_n\right), (Y_n)_{\cdot,k} = \Sigma^{-1/2}(X_n)_{\cdot,k} = \sqrt{\frac{\nu}{U_k}} Z_{\cdot,k}.$$

Noticing that $\mathbb{E}\left[\frac{1}{U_1^4}\right] = \frac{1}{(\nu-8)(\nu-6)(\nu-4)(\nu-2)}$, we have:

$$\frac{1}{p_n} \sum_{i=1}^{p_n} \mathbb{E}[y_{i1}^8] = 105 \frac{(\nu-8)(\nu-6)(\nu-4)}{(\nu-2)^3} \|\Sigma \Sigma^\top\|_n^2. \tag{2}$$

As said in the Gaussian case, in both situation we have $\|\Sigma_n \Sigma_n^\top\|_n^2 = 1$, so Assumption 2 is respected here too.

*Assumption 3.* Let $Y_n$ be $n$ iid samples of $p_n$ uncorrelated dimensions drawn from a centered Gaussian distribution, then from a Gaussian property, for all $(i, j, k, l) \in [\![1, p_n]\!]^4$ where $i, j, k, l$ are all different, we have $\text{Cov}[y_{i1}^n y_{j1}^n, y_{k1}^n y_{l1}^n] = 0$. In the case where $Y_n$ are $n$ iid samples of $p_n$ uncorrelated dimensions drawn from a centered multivariate t-distribution. Then we use the decomposition $Y = \sqrt{\frac{\nu}{U}} Z$ where $Z$ is drawn from a multivariate Gaussian distribution independent from $U$, drawn from a $\chi_\nu^2$ distribution. As for all $i \neq j$, $\text{Cov}(y_{i1}^n, y_{j1}^n) = 0$, then we trivially have $\text{Cov}(z_{i1}^n, z_{j1}^n) = 0$. So $z_{i1}$ and $z_{j1}$ are independent, which immediately leads to the fact that for all $(i, j, k, l) \in [\![1, p_n]\!]^4$ where $i, j, k, l$ are all different, we have $\text{Cov}[y_{i1}^n y_{j1}^n, y_{k1}^n y_{l1}^n] = 0$. This proves that Assumption 3 is respected in all the experimental cases we studied.

### 4.3. Results

In the following, we will use abbreviations to refer the different expected losses of each algorithms. Concerning the variants of Ledoit-Wolf shrinkage estimators with unknown mean, we denote:

- *LW_u* for the estimator we propose in this paper,
- *LW_r* for the implementation recommended by Ledoit and Wolf in 2020 [21],

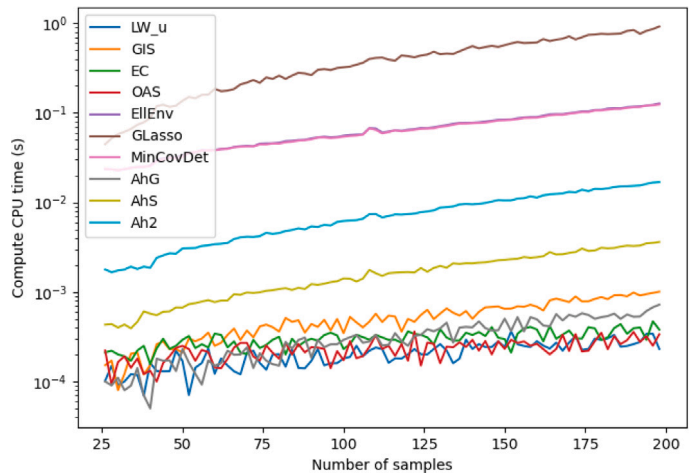

**Fig. 1.** With $p_n/n = 0.5$, random covariance, the mean compute time in seconds of each estimator in function of the number of samples $n$. All Ledoit-Wolf type estimators have very similar time complexity, only *LW_u* is shown here.

- *LW_s* for the implementation of ScikitLearn 1.2.2,
- *LW_m* for the natural estimator,
- *LW_ex* for the oracle estimator $\Sigma^*$,
- *LW_op* for the optimal estimator $\Sigma^{**}$,

Concerning the other baseline algorithms implemented in ScikitLearn, we have:

- *EC* for the Empirical Covariance estimator,
- *SC* for the Shrunk Covariance estimator,
- *OAS* for the Oracle Approximated Shrinkage estimator,

State-of-the-art baselines using non-linear shrinkage are compared too:

- *LW_an* for the Analytical non-linear shrinkage covariance estimator described in Ledoit and Wolf [17],
- *GIS* for the Geometric-Inverse Shrinkage estimator,
- *LIS* for the Linear-Inverse Shrinkage estimator,
- *QIS* for the Quadratic-Inverse Shrinkage estimator, all three described in [23].

For readability, as GIS seems to have consistently better performance among the three last estimators, we only show this one on the figures. We did not run at full scale (only with $n_{MC} = 100$) the Elliptic Envelope, GLasso, MinCovDet estimators present in ScikitLearn, and Tyler M-estimator variants as:

- *Ty* for the Tyler M-estimator with diagonal shrinkage [24],
- *SRTy* for the Tyler M-estimator with spectral shrinkage [25],
- *AhG* for the Tyler M-estimator for Gaussian distribution [13], Eq. (21),
- *AhS* for the Tyler M-estimator for Student distribution [13], Eq. (22),
- *Ah2* for the Tyler M-estimator for Student distribution [13], Eq. (31).

It is due to time complexity: in our setup, the computing time of those ones exceeds by a factor at least 10 - except AhG - the computing time of the shrinkage estimators listed before. Consequently, for reason of feasibility, we chose not to compare to them at full scale but only give an idea of how they perform, considering that the latter algorithms are part of a different class of estimators. For readability, as Ty and SRTy seems to be outperformed in our setup by AhG, AhS and Ah2, we only show the three last ones on the figures. The mean computing time in seconds are shown in Fig. 1.

### 4.3.1. Constant covariance $\Sigma_n = I_{p_n}$
*Study on a grid over (p,n).*

As they often show similar behaviors, we only show a subset of the experimental results for brevity. The three estimators LW_s, LW_r and LW_m have a very similar behavior compared to LW_u in this scenario, that is why we will only show the comparison with LW_m, having the best performance among the three. The results are shown in Fig. 2. The black contour on the surface plots is the iso-line at level 0, where the expected losses are equal. In this scenario, LW_u is constantly better than the other estimators, and the important difference is in the part $p > n$, where the mean estimation affects a lot the overall covariance estimation.

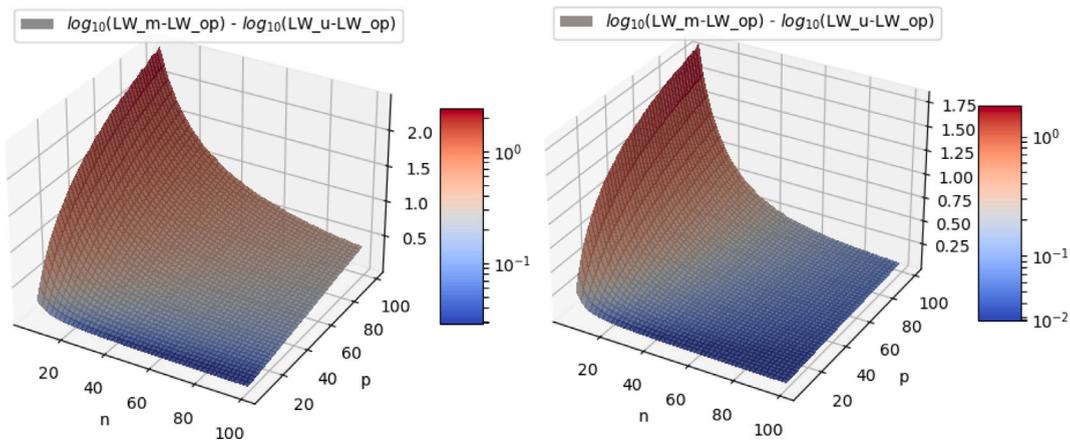

**Fig. 2.** With $p$ the dimension, $n$ the number of samples, the axis $z$ shows the difference of the $\log_{10}$ of the expected losses between LW_m and LW_u. The losses are relative to the theoretical bound LW_op. The samples are drawn with a Gaussian (left)/$t_{8.5}$ (right) distribution, $\Sigma = I_p$.

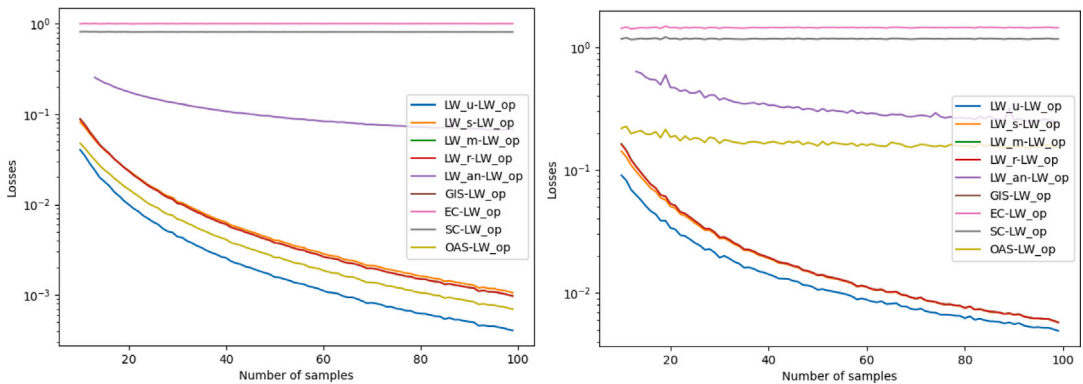

**Fig. 3.** Loss comparison to LW_op, Gaussian (left)/$t_{8.5}$ (right) distribution, $\Sigma = I_p$, $c = 1$, $n_{MC} = 10000$. Here, note that LW_ex and LW_op have both null loss, we did not plot them. Note that LW_r and LW_m are not discernible.

*Convergence study.*

We now fix $c = p_n/n$ and study the convergence of the different algorithms we cited in the experimental setup. We only show the results with the typical $c = 1$ as the other cases only widen the differences when $c$ increases and narrows then down when $c$ decreases, but do not change the order. An example with $c = 0.25$ is shown to underline this point, and show the behavior of GIS which is not defined when $c = 1$. The cases $t_{10}$ and $t_{8.5}$-distributions are very similar, that is why we show only the $t_{8.5}$ one. The key difference between the Gaussian case and the t-distribution, is that the OAS does not converge in the latter, while not being so efficient in the Gaussian case which is tailored for it. We see also a good robustness of the estimators regarding heavy-tails in the case $\nu = 4.5$. The results are shown in Figs. 3, 4, and 5.

A slight comparison with slower estimators as GLasso, MinCovDet or Tyler M-estimators is given in Fig. 6. In this setup, linear shrinkage has overall the best performance across $c$ or tails.

### 4.3.2. Covariance drawn from normalized Wishart
*Study on a grid over (p,n).*

The three estimators LW_r and LW_m have a very similar behavior compared to LW_u in this scenario, that is why we will only show the comparison with LW_m, having the best performance among the two. The results show a similar behavior between the three distributions we studied. We show the Gaussian case in Fig. 7. The black contour on the surface plots is the iso-line at level 0, it is where the losses are equal. In the case $p > n$, LW_u is far better than the other estimators, where the mean estimation affects a lot the overall covariance estimation. In a finite subset of the part $n > p$, LW_s is slightly better. LW_m presents no significant advantage compared to the two others.

*Convergence study.*

We fix $c = p_n/n$ and study the convergence of the different algorithms we cited in the experimental setup. The cases $t_{10}$ and $t_{8.5}$-distributions being very similar, we only show $t_{8.5}$. The key difference between the Gaussian case and the t-distribution, is that the OAS does not converge in the latter, while not being so efficient in the Gaussian case which is tailored for it. When $c = 1$, the differences between the linear shrinkage estimators are very low. However, key differences appear when $c$ is higher, as shown in Figs. 8, 9 with $c = 4$. A reference with $c = 0.25$ is given for comparison in Fig. 10. An interesting point is the behavior of GIS. When $c$ is low, it is the best algorithm we tested in this setup. However, it is very sensitive to heavy-tails, as $\nu = 4.5$, and to $c$ near 1, where the performances plummet. Linear shrinkage is a very robust and regularity of this behavior through different situations is an advantage for the practitioner.

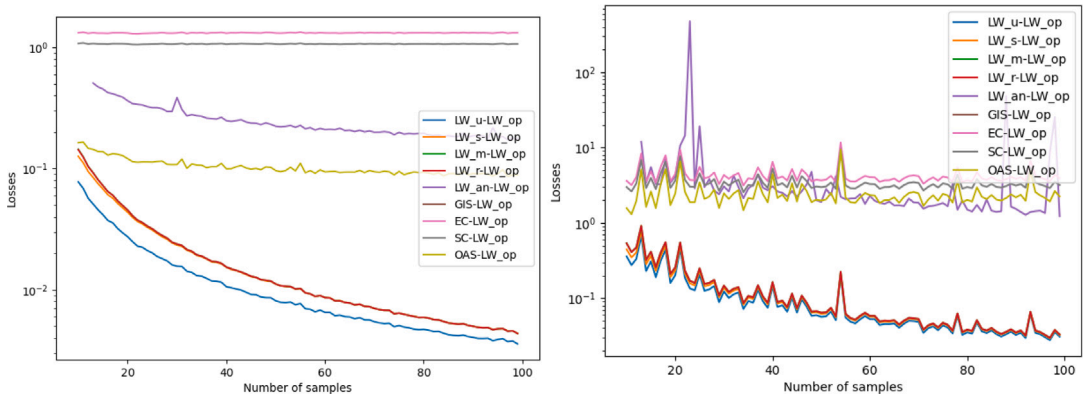

**Fig. 4.** Loss comparison to LW_op, mix of $t_{8.5}$ and $t_{15}$ (left)/$t_{4.5}$ (right) distributions, $\Sigma = I_p$, $c = 1$, $n_{MC} = 10000$. Here, note that LW_ex and LW_op have both null loss, we did not plot them. Note that LW_r and LW_m are not discernible.

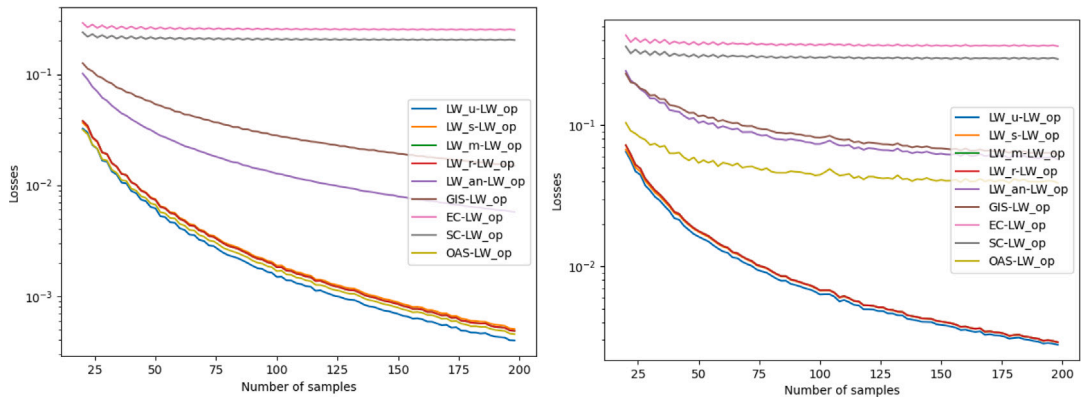

**Fig. 5.** Loss comparison to LW_op, Gaussian (left)/$t_{8.5}$ (right) distribution, $\Sigma = I_p$, $c = 0.25$, $n_{MC} = 10000$. Here, note that LW_ex and LW_op have both null loss, we did not plot them. Note that LW_r and LW_m are not discernible.

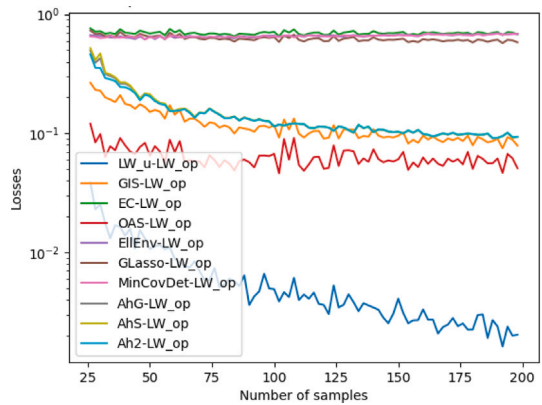

**Fig. 6.** Loss comparison to LW_op, mix of $t_{8.5}$ and $t_{15}$ distribution, $\Sigma = I_p$, $c = 0.5$, $n_{MC} = 100$.

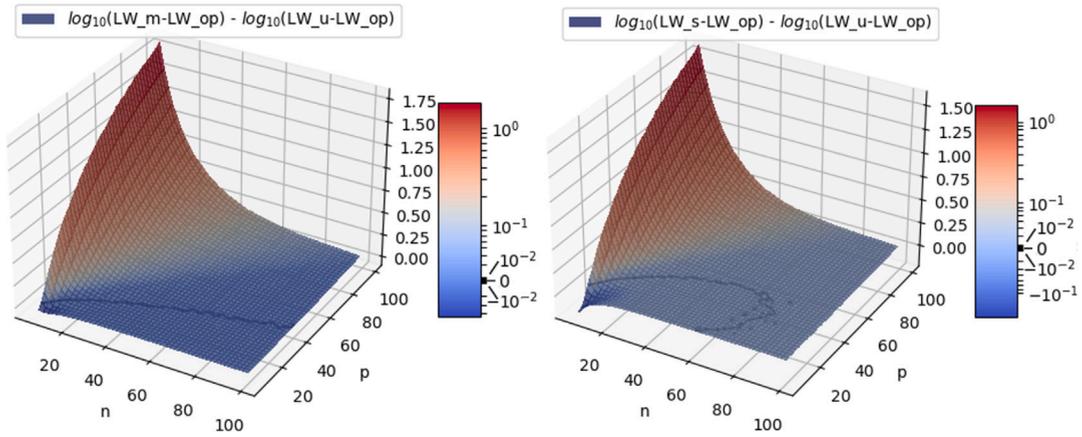

**Fig. 7.** With $p$ the dimension, $n$ the number of samples, the axis $z$ shows the difference of the $\log_{10}$ of the expected losses between LW_m (left)/LW_s (right) and LW_u. The losses are relative to the theoretical bound LW_op. Samples are drawn from a Gaussian distribution, random $\Sigma$. The black contour is the iso-line at level 0, where the expected losses are equal.

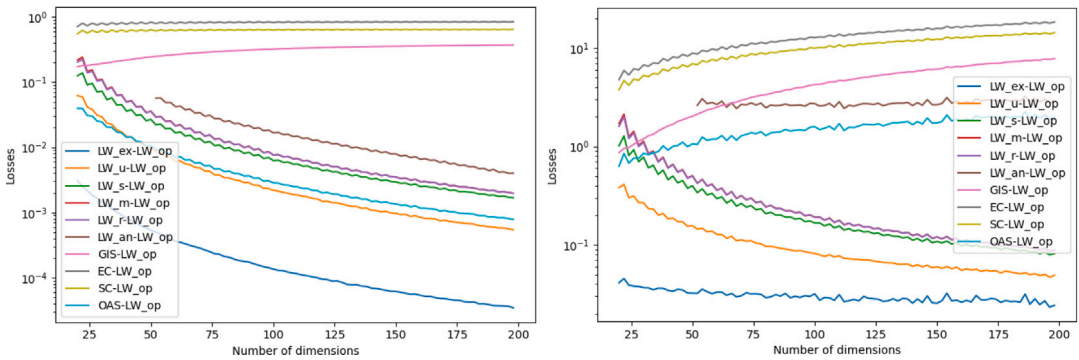

**Fig. 8.** Loss comparison to LW_op, Gaussian (left)/$t_{8.5}$ (right) distribution, random $\Sigma$, $c = 4$, $n_{MC} = 10000$.

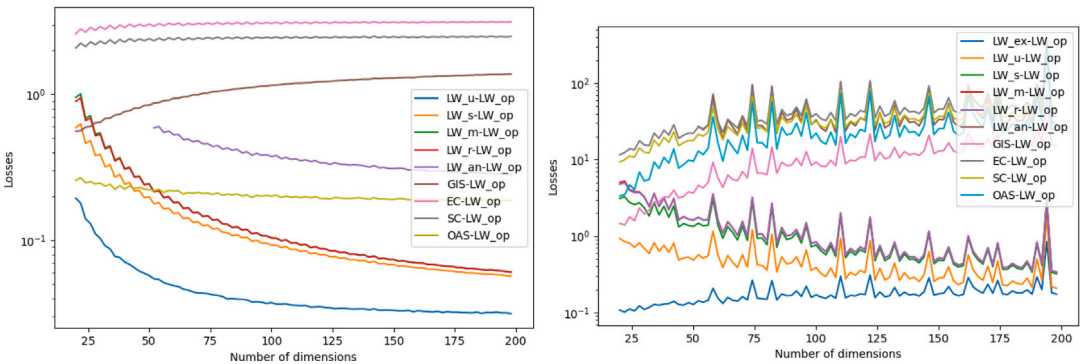

**Fig. 9.** Loss comparison to LW_op, mix of $t_{8.5}$ and $t_{15}$ (left)/$t_{4.5}$ (right) distributions, random $\Sigma$, $c = 4$, $n_{MC} = 10000$. L_ex was not computed analytically for the mixed distribution $t_{8.5}$ and $t_{15}$.

A slight comparison with slower estimators as GLasso, MinCovDet or Tyler M-estimator is given in Fig. 11. In this setup, GIS is better when $c$ is low (around 0.25 and when tails are light), however linear shrinkage has overall the best performance across $c$ or tails, underlining the robustness of this method. An interesting point is to remark that while linear shrinkage estimators are still converging towards $\Sigma^{**}$ when $\nu = 4.5$, it is not the case of $\Sigma^*$ anymore, which has a significantly worse loss. Our estimators, having random coefficients, outperform the non-random oracle in this heavy-tail situation.

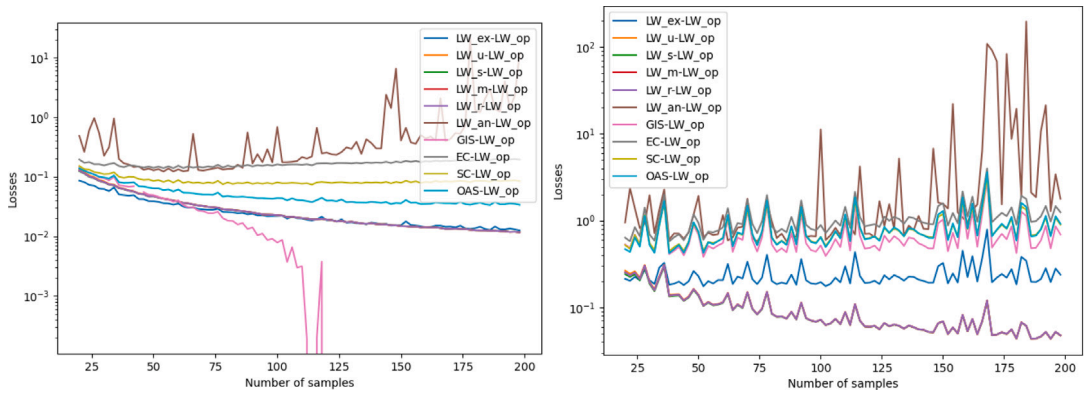

**Fig. 10.** Loss comparison to LW_op, $t_{8.5}$ (left)/$t_{4.5}$ (right) distribution, random $\Sigma$, $c = 0.25$, $n_{MC} = 10000$. GIS slightly outperforms LW_op in the left situation and in the Gaussian case, but is sensible to heavy tails on the right.

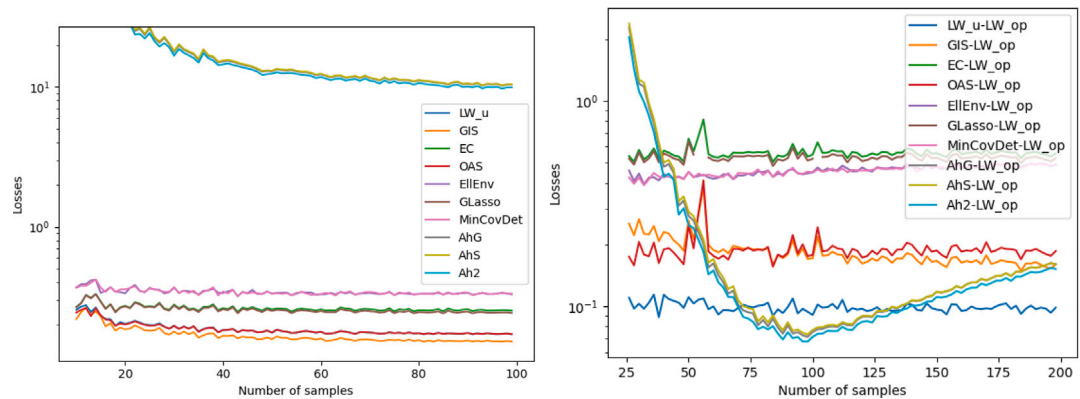

**Fig. 11.** Loss comparison to LW_op, Gaussian (left)/ mix of $t_{8.5}$ and $t_{15}$ (right) distribution, random $\Sigma$, $c = 0.25$ (left)/$c = 0.5$ (right), $n_{MC} = 100$.

## 5. Conclusion

In this work, we extended the linear shrinkage approach of Ledoit and Wolf [10] for covariance matrix estimation to the case where the mean of the distribution is unknown. Theoretically, we showed that in this case we have similar asymptotic properties as in the situation when the mean is known. Four different translation-invariant estimators emerged, three around those implemented in ScikitLearn or recommended by Ledoit and Wolf, and one naturally emerging from the theoretical proofs. Experimentally, the latter showed improved performances in a large spectrum of situations compared to ScikitLearn 1.2.2 baselines, non-linear shrinkage estimators, M-estimators, and to the three other estimators presented in the theoretical part. The gain in performance is particularly high when the dimension is bigger than the number of samples, while the differences are comparably low when dimension is smaller than the number of samples.

Linear shrinkage methods in this setup of unknown mean shows a robustness to different types of distributions, tails and concentrations, and remains state-of-the-art in his results facing a variety of situations. The theoretical extensions and experimental conclusions we had in this work aims at giving the practitioner keys to choose the right estimator regarding its constraints.

Ledoit and Wolf developed several non-linear shrinkage estimators where the mean is known [15–17]. Work needs to be conducted to investigate if similar approach can be used to extend their non-linear frameworks.

**CRediT authorship contribution statement**

**Benoît Oriol:** Conceptualization, Formal analysis, Investigation, Methodology, Visualization Writing – original draft, Writing – review & editing. **Alexandre Miot:** Formal analysis, Writing – review & editing.

**Acknowledgments**

We thank the Editor, Associate Editor and referees, Gabriel Turinici for its advices all along the work as well as the ANRT for its funding through the CIFRE contract 2022/0536.

**Appendix**

In this Appendix, we present the proofs of the technical results. For brevity, we omit the subscript $n$; but it is understood that everything depends on $n$. Coefficients of $\Lambda$ are denoted $\lambda_{ij}$, and if not stated otherwise, sum indices $i, j, l, m$ are in $[\![1, p]\!]$ and $k_1, k'_1, k_2, k'_2, \ldots$ are in $[\![1, n]\!]$. Moreover, we denote $\bar{y} = \frac{1}{n} \sum_k y_{\cdot,k}$. We recall that, from Remark 2, as all the estimators are invariant by change of mean from the definition of $\tilde{X} = X - \sum_k X_{\cdot,k}$, we assume $\mathbb{E}[X] = 0$ for the simplicity of notation.

*A.1. Proof of Lemma 1*

We have:

$$\|\Sigma\|^2 = \|\Lambda\|^2 = \frac{1}{p} \sum_{i=1}^{p} \mathbb{E}[y_{i1}^2]^2 \leq \frac{1}{p} \sum_{i=1}^{p} \mathbb{E}[y_{i1}^4] \leq \sqrt{\frac{1}{p} \sum_{i=1}^{p} \mathbb{E}[y_{i1}^8]} \leq \sqrt{K_2}. \tag{3}$$

As $\mu = \langle \Sigma, I \rangle \leq \|\Sigma\|$, $\mu$ remains bounded as $n$ goes to infinity.

Also, $\alpha^2 = \|\Sigma - \mu I\|^2 = \|\Sigma\|^2 - \mu^2$, so remains bounded as $n$ goes to infinity too.

For $\beta^2$, we will deeply decompose the expectation. This is not absolutely necessary to prove the boundedness, but the decomposition will be of utter importance in the following proofs. So, we have by reindexation:

$$\beta^2 = \frac{1}{p} \sum_{i,j} \frac{1}{n^2(n-1)^2} \sum_{k_1} \sum_{k'_1 \neq k_1} \sum_{k_2} \sum_{k'_2 \neq k_2} \mathbb{E}\left[ \left( y_{ik_1}(y_{jk_1} - y_{jk'_1}) - \lambda_{ij} \right) \left( y_{ik_2}(y_{jk_2} - y_{jk'_2}) - \lambda_{ij} \right) \right]. \tag{4}$$

We denote, for $k_1, k'_1 \neq k_1, k_2, k'_2 \neq k_2$:

$$E_{ij}(k_1, k'_1, k_2, k'_2) = \mathbb{E}\left[ \left( y_{ik_1}(y_{jk_1} - y_{jk'_1}) - \lambda_{ij} \right) \left( y_{ik_2}(y_{jk_2} - y_{jk'_2}) - \lambda_{ij} \right) \right]. \tag{5}$$

- If $|\{k_1, k'_1\} \cap \{k_2, k'_2\}| = 0$: $E_{ij}(k_1, k'_1, k_2, k'_2) = 0$.
- If $|\{k_1, k'_1\} \cap \{k_2, k'_2\}| = 1$:

  – If $k_1 = k_2$: then $k'_1 \neq k'_2$, and $k_1 \neq k'_2$. So, $E_{ij}(k_1, k'_1, k_1, k'_2) = \mathbb{V}[y_{i1} y_{j1}]$. There are $n(n-1)(n-2)$ terms respecting those conditions.
  – Else: $E_{ij}(k_1, k'_1, k_2, k'_1) = 0$.

- If $|\{k_1, k'_1\} \cap \{k_2, k'_2\}| = 2$:

  – If $k_1 = k_2$ and $k'_1 = k'_2$: then $k_1 \neq k'_1$. So, $E_{ij}(k_1, k'_1, k_1, k'_1) = \mathbb{V}[y_{i1} y_{j1}] + \lambda_{ii} \lambda_{jj}$. There are $n(n-1)$ terms respecting those conditions.
  – If $k_1 = k'_2$ and $k'_1 = k_2$: then $k_1 \neq k'_1$. So, $E_{ij}(k_1, k'_1, k'_1, k_1) = \lambda_{ij}^2$. There are $n(n-1)$ terms respecting those conditions.

Using the latter decomposition on $(k_1, k'_1, k_2, k'_2)$, we deduce:

$$\beta^2 = \frac{1}{pn} \sum_{i,j} \mathbb{V}[y_{i1} y_{j1}] + \frac{p+1}{n(n-1)} \mu^2 + \frac{1}{n(n-1)} \alpha^2. \tag{6}$$

$\mu^2$ and $\alpha^2$ are bounded and following Ledoit and Wolf proof of Lemma 3.1 [10],

$\frac{1}{pn} \sum_{i,j} \mathbb{V}[y_{i1} y_{j1}] \leq K_1 \sqrt{K_2}$, so $\beta^2$ remains bounded as $n$ goes to infinity. Finally, $\delta^2 = \alpha^2 + \beta^2$ is also bounded as $n$ goes to infinity, which conclude the proof of the lemma.

*A.2. Proof of Theorem 1*

Reindexing, and using the proof of Lemma 1 for the upper bound, we have:

$$\mu^2 + \theta^2 = \frac{1}{p}(\alpha^2 + \mu^2) + \frac{1}{p^2} \sum_{ij} \mathbb{V}\left[ y_{i1} y_{j1} \right] \leq \frac{1}{p}(\alpha^2 + \mu^2) + \sqrt{K_2}. \tag{7}$$

As $\alpha^2$ and $\mu^2$ remains bounded as $n$ goes to infinity, so is $\theta^2$. And, from the proof of Lemma 1:

$$\beta^2 = \frac{1}{pn} \sum_{i,j} \mathbb{V}[y_{i1} y_{j1}] + \frac{p+1}{n(n-1)} \mu^2 + \frac{1}{n(n-1)} \alpha^2 = \frac{p}{n}(\mu^2 + \theta^2) - \frac{1}{n(n-1)} \alpha^2 + \frac{p-n+2}{n(n-1)} \mu^2. \tag{8}$$

$\alpha^2$, $\mu^2$ and $\frac{p}{n}$ remains bounded as $n$ goes to infinity, which conclude the proof of the theorem.

## A.3. Proof of Lemma 2

By linearity of the inner product, we trivially have: $\mathbb{E}[m] = \mu$. For the quartic mean convergence, we decompose and reindex:

$$\mathbb{E}[(m-\mu)^4] = \frac{1}{(n-1)^4} \sum_{\substack{k_1,k_2,k_3,k_4 \\ k'_1,k'_2,k'_3,k'_4}} \mathbb{E}\left[\prod_{s=1}^{4}\left(\frac{1}{pn}\sum_i (y_{ik_s}^2 - \lambda_{ii})\right)\right]$$

$$- 4\mathbb{E}\left[\prod_{s=1}^{3}\left(\frac{1}{pn}\sum_i (y_{ik_s}^2 - \lambda_{ii})\right)\left(\frac{1}{pn}\sum_i (y_{ik_4}y_{ik'_4} - \delta_{k_4=k'_4}\lambda_{ii})\right)\right]$$

$$+ 6\mathbb{E}\left[\prod_{s=1}^{2}\left(\frac{1}{pn}\sum_i (y_{ik_s}^2 - \lambda_{ii})\right)\prod_{t=3}^{4}\left(\frac{1}{pn}\sum_i (y_{ik_t}y_{ik'_t} - \delta_{k_t=k'_t}\lambda_{ii})\right)\right] \quad (9)$$

$$- 4\mathbb{E}\left[\left(\frac{1}{pn}\sum_i (y_{ik_1}^2 - \lambda_{ii})\right)\prod_{t=2}^{4}\left(\frac{1}{pn}\sum_i (y_{ik_t}y_{ik'_t} - \delta_{k_t=k'_t}\lambda_{ii})\right)\right]$$

$$+ \mathbb{E}\left[\prod_{t=1}^{4}\left(\frac{1}{pn}\sum_i (y_{ik_t}y_{ik'_t} - \delta_{k_t=k'_t}\lambda_{ii})\right)\right].$$

If $k_1, k'_1, k_2, k'_2, k_3, k'_3, k_4, k'_4$ are all different, then the expectation in the sum equals 0. So there are at most $28n^7 + O(n^6)$ non-zero terms in the sum.

Now, let us find a bound of those expectations. Let us first note that, for all $(k, k') \in [\![1, n]\!]^2$ and $i \in [\![1, p]\!]$:

$$\mathbb{E}\left[\left(y_{ik}y_{ik'} - \delta_{k=k'}\lambda_{ii}\right)^4\right] = \mathbb{E}[(y_{ik}y_{ik'})^4] - 4\mathbb{E}[(y_{ik}y_{ik'})^3]\mathbb{E}[y_{ik}y_{ik'}] + 6\mathbb{E}[(y_{ik}y_{ik'})^2]\mathbb{E}[y_{ik}y_{ik'}]^2 - 3\mathbb{E}[y_{ik}y_{ik'}]^4. \quad (10)$$

If $k = k'$ then $\mathbb{E}[(y_{ik}y_{ik'})^3]\mathbb{E}[y_{ik}y_{ik'}] \geq 0$ and if $k \neq k'$, $\mathbb{E}[(y_{ik}y_{ik'})^3]\mathbb{E}[y_{ik}y_{ik'}] = 0 \geq 0$, so by Cauchy–Schwarz: $\mathbb{E}\left[\left(y_{ik}y_{ik'} - \delta_{k=k'}\lambda_{ii}\right)^4\right] \leq 7\mathbb{E}[(y_{ik}y_{ik'})^4]$. Back to the bound of our expectation, let $N \in [\![1, 4]\!]$, and using multiple Cauchy–Schwarz and Jensen inequalities:

$$\left|\mathbb{E}\left[\prod_{s=1}^{N}\left(\frac{1}{pn}\sum_i (y_{ik_s}^2 - \lambda_{ii})\right)\prod_{t=N+1}^{4}\left(\frac{1}{pn}\sum_i (y_{ik_t}y_{ik'_t} - \delta_{k_t=k'_t}\lambda_{ii})\right)\right]\right| \leq \frac{7K_2}{n^4}. \quad (11)$$

So, in conclusion of this proof,

$$\mathbb{E}[(m-\mu)^4] \leq 7(1 + 4 + 6 + 4 + 1)\frac{28n^7 + O(n^6)}{n^4(n-1)^4} K_2 \xrightarrow[n \to \infty]{} 0. \quad (12)$$

## A.4. Proof of Corollary 2

We have:

$$\mathbb{V}[m^2] = \mathbb{E}[(m-\mu)^4] + 4\mathbb{E}[(m-\mu)^3]\mu + 4\mathbb{E}[(m-\mu)^2]\mu^2 - \mathbb{E}[(m-\mu)^2]^2. \quad (13)$$

And, from Lemma 2, $\mathbb{E}[(m-\mu)^4] \to 0$, and so goes for the smaller moments by Jensen inequality, and $\mu$ is bounded from Lemma 1. So, $\mathbb{V}[m^2] \to 0$.

## A.5. Proof of Lemma 3

### A.5.1. Preliminary combinatorial result

Let $K \in \mathbb{N}^*$, and $K$ indices $(k_1, \ldots, k_K) \in [\![1, n]\!]^K$.

Let us associate a graph with $K$ vertices $\mathcal{V} = \{1, \ldots, K\}$ to this set of indices. The set of edges $\mathcal{E}$ is built as following: there is an edge between the node $a \in \mathcal{V}$ and $b \in \mathcal{V}$, $a \neq b$ (we do not allow self-loops), if the corresponding indices are equal, i.e if $k_a = k_b$. We finally define our graph $\mathcal{G} = (\mathcal{V}, \mathcal{E})$.

**Proposition 1.** *Let $\mathcal{G} = (\mathcal{V}, \mathcal{E})$ a graph with $K$ vertices generated from some indices $(k_1^{(0)}, \ldots, k_K^{(0)}) \in [\![1, n]\!]^K$ with the procedure described previously. Suppose $\mathcal{G}$ has $C \in [\![1, K]\!]$ connected components. Then, there are $\prod_{i=0}^{C-1}(n-i)$ set of indices $(k_1, \ldots, k_K) \in [\![1, n]\!]^K$ which have the associated graph $\mathcal{G}$.*

For each node $v \in \mathcal{V}$, $v$ belongs to a unique connected component that we denote $c(v) \in [\![1, C]\!]$. Then, the function:

$$x \in \left\{x \in [\![1, n]\!]^C, x_1, \ldots, x_C \text{ all different}\right\} \mapsto (x_{c(1)}, \ldots, x_{c(K)})$$

is a bijection between $\left\{x \in [\![1, n]\!]^C, x_1, \ldots, x_C \text{ all different}\right\}$ and $\left\{(k_1, \ldots, k_K) \in [\![1, n]\!]^K \text{ which have the associated graph } \mathcal{G}\right\}$. Immediately, we deduce that its cardinal is equal to $\prod_{i=0}^{C-1}(n-i)$.

*A.5.2. Proof of Lemma 3*

From the proof of Lemma 3.3 in Ledoit and Wolf [10], we have:

$$d^2 - \delta^2 = -(m-\mu)^2 + 2\mu(\mu - m) + \left(\|S\|^2 - \mathbb{E}[\|S\|^2]\right). \tag{14}$$

The first two terms converge to 0 in quadratic mean thanks to Lemma 2. Let us show that the last term converges to 0 in quadratic mean too, i.e $\mathbb{V}[\|S\|^2] \xrightarrow[n\to\infty]{} 0$.

Decomposing $\|S\|^2$, we have:

$$\begin{aligned}
\|S\|^2 &= \frac{p}{n^2} \sum_{k_1} \left(\frac{1}{p}\sum_i y_{ik_1}^2\right)^2 & (t1') \\
&+ \frac{p(n^2 - 2n + 2)}{n^2(n-1)^2} \sum_{k_1, k_1' \neq k_1} \left(\frac{1}{p}\sum_i y_{ik_1} y_{ik_1'}\right)^2 & (t2') \\
&- \frac{p(2n-3)}{n^2(n-1)^2} \sum_{k_1, k_1' \neq k_1, k_2' \neq k_1'} \frac{1}{p^2}\sum_{i,j} y_{ik_1} y_{ik_1'} y_{jk_1} y_{jk_2'} & (t3') \\
&- \frac{p(2n-1)}{n^2(n-1)^2} \sum_{k_1, k_1' \neq k_1} \frac{1}{p^2}\sum_{i,j} y_{ik_1} y_{ik_1'} y_{jk_1'}^2 & (t4') \\
&+ \frac{p}{n^2(n-1)^2} \sum_{\substack{k_1, k_1' \neq k_1 \\ k_2 \neq k_1', k_2' \neq k_2}} \frac{1}{p^2}\sum_{i,j} y_{ik_1} y_{ik_2} y_{jk_1'} y_{jk_2'} & (t5').
\end{aligned} \tag{15}$$

It is sufficient to show that the variance of each of the 5 term converges to 0 as $n$ goes to infinity in order to prove that $\mathbb{V}[\|S\|^2]$ converges to 0 as $n$ goes to infinity.

From the proof of Lemma 3.3 in (Ledoit and Wolf, 2004) [10], we immediately have that $\mathbb{V}[(t1')]$ and $\mathbb{V}[(t2')]$ converge to 0.

For the three remaining terms, we detail only the more complex one, i.e., $\mathbb{V}[(t3')]$, as the method apply directly the other terms.

$$\mathbb{V}[(t3')] = \frac{p^2(2n-3)^2}{n^4(n-1)^4} \sum_{\substack{k_1, k_1' \neq k_1 \\ k_2' \neq k_1'}} \sum_{\substack{k_3, k_3' \neq k_3 \\ k_4' \neq k_3'}} \sum_{i,j,\ell,m} \frac{1}{p^4} \mathrm{Cov}(y_{ik_1} y_{ik_1'} y_{jk_1} y_{jk_2'}, y_{\ell k_3} y_{\ell k_3'} y_{mk_3} y_{mk_4'}). \tag{16}$$

Let $(k_1, k_1', k_2', k_3, k_3', k_4') \in [\![1,n]\!]^6$ respecting the conditions given in the sums. Suppose there exists $(i,j,\ell,m) \in [\![1,p]\!]^4$ such that $\mathrm{Cov}(y_{ik_1} y_{ik_1'} y_{jk_1} y_{jk_2'}, y_{\ell k_3} y_{\ell k_3'} y_{mk_3} y_{mk_4'}) \neq 0$. Let consider the graph $\mathcal{G} = (\mathcal{V}, \mathcal{E})$ built from $(k_1, k_1', k_2', k_3, k_3', k_4')$ following the procedure described in the preliminary combinatorial result.

As for some $(i,j,\ell,m) \in [\![1,p]\!]^4$ we have $\mathrm{Cov}(y_{ik_1} y_{ik_1'} y_{jk_1} y_{jk_2'}, y_{\ell k_3} y_{\ell k_3'} y_{mk_3} y_{mk_4'}) \neq 0$, by independence, the nodes 1', 2', 3', and 4' cannot be isolated. As a consequence, $\mathcal{G}$ has at least 2 edges.

- When the graph $\mathcal{G}$ has only 2 edges, we have either of the following conditions, that we denote the (*) conditions:

  - $(k_1' = k_3') \wedge (k_2' = k_4') \wedge (k_1 \neq k_3) \wedge (k_1 \neq k_2') \wedge (k_3 \neq k_4')$, or
  - $(k_1' = k_4') \wedge (k_3' = k_2') \wedge (k_1 \neq k_3) \wedge (k_1 \neq k_2') \wedge (k_3 \neq k_4')$.

Note that the case where $(k_1' = k_2') \wedge (k_3' = k_4')$ is impossible due to the constraints on the indices in the sum. In the first case, we have:

$$\mathrm{Cov}(y_{ik_1} y_{ik_1'} y_{jk_1} y_{jk_2'}, y_{\ell k_3} y_{\ell k_3'} y_{mk_3} y_{mk_4'}) = \lambda_{ij}\lambda_{\ell m}\lambda_{i\ell}\lambda_{jm}. \tag{17}$$

And in the second case, we have:

$$\mathrm{Cov}(y_{ik_1} y_{ik_1'} y_{jk_1} y_{jk_2'}, y_{\ell k_3} y_{\ell k_3'} y_{mk_3} y_{mk_4'}) = \lambda_{ij}\lambda_{\ell m}\lambda_{im}\lambda_{j\ell}. \tag{18}$$

Using the fact that $i \neq j \implies \lambda_{ij} = 0$, in both cases we have:

$$\sum_{i,j,l,m} \mathrm{Cov}(y_{ik_1} y_{ik_1'} y_{jk_1} y_{jk_2'}, y_{\ell k_3} y_{\ell k_3'} y_{mk_3} y_{mk_4'}) \leq pK_2. \tag{19}$$

Under the conditions (*), $\mathcal{G}$ has exactly 4 connected components. So, from the preliminary combinatorial result, there are $2n(n-1)(n-2)(n-3)$ different 6-uples of indices $k$ respecting the (*) conditions. We finally have:

$$\begin{aligned}
&\frac{p^2(2n-3)^2}{n^4(n-1)^4} \sum_{\substack{(k_1,k_1',k_2',k_3,k_3',k_4') \\ \text{under (*) conditions}}} \sum_{i,j,\ell,m} \frac{1}{p^4} \mathrm{Cov}(y_{ik_1} y_{ik_1'} y_{jk_1} y_{jk_2'}, y_{\ell k_3} y_{\ell k_3'} y_{mk_3} y_{mk_4'}) \\
&\leq \frac{p^2(2n-3)^2}{n^4(n-1)^4} \times 2n(n-1)(n-2)(n-3) \frac{1}{p^4} pK_2 \xrightarrow[n\to\infty]{} 0.
\end{aligned} \tag{20}$$

- Otherwise, $\mathcal{G}$ has 3 edges or more: we denote it as the (**) condition. As it has only 6 vertices, there are at most 3 connected component. So, from the preliminary combinatorial result, there are $n(n-1)(n-2) + n(n-1) + n = n((n-1)^2 + n)$ different 6-uples of indices $k$ that have $\mathcal{G}$ as associated graph. Moreover, there are a finite number $N \in \mathbb{N}$ (independent of $n$) of graphs with 6 vertices that have at most 3 connected components.

So, there are at most $Nn((n-1)^2 + n)$ different 6-uples of indices $k$ such that the associated graph $\mathcal{G}$ has 3 edges or more. And we have, by Cauchy–Schwarz and Jensen inequalities:

$$\sum_{i,j,\ell,m} \frac{1}{p^4} \left| \text{Cov}(y_{ik_1} y_{ik'_1} y_{jk_1} y_{jk'_2}, y_{\ell k_3} y_{\ell k'_3} y_{mk_3} y_{mk'_4}) \right| \leq K_2. \tag{21}$$

So, using both of the previous inequalities,

$$\frac{p^2(2n-3)^2}{n^4(n-1)^4} \sum_{\substack{(k_1, k'_1, k_2, k_3, k'_3, k'_4) \\ \text{under (**) condition}}} \sum_{i,j,\ell,m} \frac{1}{p^4} \text{Cov}(y_{ik_1} y_{ik'_1} y_{jk_1} y_{jk'_2}, y_{\ell k_3} y_{\ell k'_3} y_{mk_3} y_{mk'_4})$$

$$\leq \frac{(2n-3)^2}{n^2(n-1)^4} \times Nn\left((n-1)^2 + n\right) K_1^2 K_2 \xrightarrow[n \to \infty]{} 0. \tag{22}$$

So, from both previous cases, we immediately have $\mathbb{V}[(t3')] \xrightarrow[n \to \infty]{} 0$. With the same method, we have: $\mathbb{V}[(t4')] \xrightarrow[n \to \infty]{} 0$ and $\mathbb{V}[(t5')] \xrightarrow[n \to \infty]{} 0$.

We showed that each of the 5 terms of $\|S\|^2$ have a variance that converges to 0 as $n$ goes to infinity.

So, $\mathbb{V}[\|S\|^2] \xrightarrow[n \to \infty]{} 0$, which concludes the proof of the first part of the lemma: $d^2 - \delta^2 \xrightarrow[n \to \infty]{q.m} 0$. Finally, by property of the expectation, we have that $\mathbb{E}[(d^2 - \mathbb{E}[d^2])^2] \leq \mathbb{E}[(d^2 - \delta^2)^2]$, so it follows that: $d^2 - \mathbb{E}[d^2] \xrightarrow[n \to \infty]{q.m} 0$.

## A.6. Proof of Lemma 4

Developing $\mathbb{E}[\bar{b}^2]$, we have:

$$\begin{aligned}
\mathbb{E}[\bar{b}^2] = \frac{n}{p(n-1)^2} \sum_{i,j} \Bigg( & \mathbb{V}\left[y_{i1} y_{j1}\right] + \mathbb{V}\left[\bar{y}_i y_{j1}\right] + \mathbb{V}\left[\bar{y}_j y_{i1}\right] + \mathbb{V}\left[\bar{y}_i \bar{y}_j\right] \\
& - 2\mathbb{E}\left[(y_{i1} y_{j1} - \lambda_{ij})\left(\bar{y}_i y_{j1} - \frac{1}{n}\lambda_{ij}\right)\right] - 2\mathbb{E}\left[(y_{i1} y_{j1} - \lambda_{ij})\left(\bar{y}_j y_{i1} - \frac{1}{n}\lambda_{ij}\right)\right] \\
& + 2\mathbb{E}\left[(y_{i1} y_{j1} - \lambda_{ij})\left(\bar{y}_i \bar{y}_j - \frac{1}{n}\lambda_{ij}\right)\right] + 2\mathbb{E}\left[\left(\bar{y}_i y_{j1} - \frac{1}{n}\lambda_{ij}\right)\left(\bar{y}_j y_{i1} - \frac{1}{n}\lambda_{ij}\right)\right] \\
& - 2\mathbb{E}\left[\left(\bar{y}_i y_{j1} - \frac{1}{n}\lambda_{ij}\right)\left(\bar{y}_i \bar{y}_j - \frac{1}{n}\lambda_{ij}\right)\right] - 2\mathbb{E}\left[\left(\bar{y}_j y_{i1} - \frac{1}{n}\lambda_{ij}\right)\left(\bar{y}_i \bar{y}_j - \frac{1}{n}\lambda_{ij}\right)\right] \Bigg) - \frac{1}{n}\beta^2.
\end{aligned} \tag{23}$$

Computing each term in function of $\mathbb{V}[y_{i1} y_{j1}]$, $\lambda_{ij}$, $\lambda_{ii}$ and $\lambda_{jj}$, we obtain:

$$\mathbb{E}\left[\bar{b}^2\right] = \frac{1}{p(n-1)} \left( \frac{n^2 - 3n + 3}{n^2} \sum_{i,j} \mathbb{V}[y_{i1} y_{j1}] + p^2 \frac{2n-3}{n^2} \mu^2 + p \frac{2n-3}{n^2} (\alpha^2 + \mu^2) \right) - \frac{1}{n}\beta^2. \tag{24}$$

From the proof of Lemma 1, we have: $\frac{1}{pn} \sum_{i,j} \mathbb{V}[y_{i1} y_{j1}] = \beta^2 - \frac{p+1}{n(n-1)} \mu^2 - \frac{1}{n(n-1)} \alpha^2$. Denoting:

$\gamma_n = \frac{n(n-1)}{n^2 - 3n + 3}$, $\lambda_n = \frac{n^2(n-2)}{(n-1)(n^2 - 3n + 3)}$, $c_0 = \frac{1}{\gamma_n} - \frac{1}{n} - \frac{\lambda_n}{\gamma_n n^2}$, $c_1 = \frac{\lambda_n}{\gamma_n n^2}$, $c_2 = (p+1)c_1$, we obtain: $\mathbb{E}[\bar{b}^2] = (c_0 + c_1)\beta^2 + c_1 \alpha^2 + c_2 \mu^2$.

So, we can conclude using $\delta^2 = \alpha^2 + \beta^2$ that $\mathbb{E}[\bar{b}^2] = c_0 \beta^2 + c_1 \delta^2 + c_2 \mu^2$.

## A.7. Proof of Lemma 5

We compute the expectations of $m^2$ and $d^2$: $\mathbb{E}[m^2] = \mu^2 + \mathbb{V}[m]$, and $\mathbb{E}[d^2] = \delta^2 - \mathbb{V}[m]$. Moreover, from Lemma 4, we have: $\mathbb{E}[\bar{b}^2] = c_0 \beta_n^2 + c_1 \delta_n^2 + c_2 \mu_n^2$. So, combining the last 3 equations, we obtain: $\mathbb{E}[\bar{b}^2] = c_0 \beta^2 + c_1 \mathbb{E}[d^2] + c_2 \mathbb{E}[m^2] + (c_1 - c_2)\mathbb{V}[m]$.

## A.8. Proof of Lemma 6

$$\bar{b}^2 - \mathbb{E}[\bar{b}^2] = \left[ \left( \frac{1}{n^2} \sum_k \left\| \frac{n}{n-1} \tilde{x}_{\cdot k} \tilde{x}_{\cdot k}^\top - \Sigma \right\|^2 \right) - \mathbb{E}[\bar{b}^2] \right] + \left[ \left( \frac{1}{n^2} \sum_k \left\| \frac{n}{n-1} \tilde{x}_{\cdot k} \tilde{x}_{\cdot k}^\top - S \right\|^2 \right) - \left( \frac{1}{n^2} \sum_k \left\| \frac{n}{n-1} \tilde{x}_{\cdot k} \tilde{x}_{\cdot k}^\top - \Sigma \right\|^2 \right) \right]. \tag{25}$$

Firstly, we want to show that the variance of the first term converges to 0 a $n$ goes to infinity. Following the decomposition developed in Lemma 4, we have:

$$\mathbb{V}\left[\frac{1}{n^2}\sum_k\left\|\frac{n}{n-1}\tilde{x}_{\cdot k}\tilde{x}_{\cdot k}^\top - \Sigma\right\|^2\right] = \mathbb{V}\left[\frac{1}{p(n-1)^2}\sum_k\sum_{i,j}\left([y_{ik}y_{jk} - \lambda_{ij}]^2 + \left[\bar{y}_i y_{jk} - \frac{\lambda_{ij}}{n}\right]^2 + \left[\bar{y}_j y_{ik} - \frac{\lambda_{ij}}{n}\right]^2 + \left[\bar{y}_i \bar{y}_j - \frac{\lambda_{ij}}{n}\right]^2\right.\right.$$
$$- 2\left[(y_{ik}y_{jk} - \lambda_{ij})\left(\bar{y}_i y_{jk} - \frac{\lambda_{ij}}{n}\right)\right] - 2\left[(y_{ik}y_{jk} - \lambda_{ij})\left(\bar{y}_j y_{ik} - \frac{\lambda_{ij}}{n}\right)\right]$$
$$+ 2\left[(y_{ik}y_{jk} - \lambda_{ij})\left(\bar{y}_i \bar{y}_j - \frac{\lambda_{ij}}{n}\right)\right] + 2\left[\left(\bar{y}_i y_{jk} - \frac{\lambda_{ij}}{n}\right)\left(\bar{y}_j y_{ik} - \frac{\lambda_{ij}}{n}\right)\right]$$
$$\left.\left.- 2\left[\left(\bar{y}_i y_{jk} - \frac{\lambda_{ij}}{n}\right)\left(\bar{y}_i \bar{y}_j - \frac{\lambda_{ij}}{n}\right)\right] - 2\left[\left(\bar{y}_j y_{ik} - \frac{\lambda_{ij}}{n}\right)\left(\bar{y}_i \bar{y}_j - \frac{\lambda_{ij}}{n}\right)\right]\right)\right].$$
(26)

We then prove that the variance of each of the 10 separated sums converges to 0. For that, let us show a useful inequality.

*A.8.1. Preliminary inequality*

Let $(k_1, k_1', k_2, k_2', k_3, k_3', k_4, k_4') \in [\![1, n]\!]^8$ and $(i, j, \ell, m) \in [\![1, p]\!]^4$. Then, chaining Cauchy–Schwarz and Jensen inequalities,

$$\mathrm{Cov}\left(\left(y_{ik_1}y_{jk_1'} - \lambda_{ij}\mathbf{1}_{k_1=k_1'}\right)\left(y_{ik_2}y_{jk_2'} - \lambda_{ij}\mathbf{1}_{k_2=k_2'}\right), \left(y_{\ell k_3}y_{mk_3'} - \lambda_{\ell m}\mathbf{1}_{k_3=k_3'}\right)\left(y_{\ell k_4}y_{mk_4'} - \lambda_{\ell m}\mathbf{1}_{k_4=k_4'}\right)\right)$$
$$\leq 9\sqrt[4]{\mathbb{E}[y_{i1}^8]\mathbb{E}[y_{j1}^8]\mathbb{E}[y_{\ell 1}^8]\mathbb{E}[y_{m1}^8]}.$$
(27)

So, for all $(k_1, k_1', k_2, k_2', k_3, k_3', k_4, k_4') \in [\![1, n]\!]^8$, we have,

$$\sum_{i,j,\ell,m}\frac{1}{p^4}\mathrm{Cov}\left(\left(y_{ik_1}y_{jk_1'} - \lambda_{ij}\mathbf{1}_{k_1=k_1'}\right)\left(y_{ik_2}y_{jk_2'} - \lambda_{ij}\mathbf{1}_{k_2=k_2'}\right), \left(y_{\ell k_3}y_{mk_3'} - \lambda_{\ell m}\mathbf{1}_{k_3=k_3'}\right)\left(y_{\ell k_4}y_{mk_4'} - \lambda_{\ell m}\mathbf{1}_{k_4=k_4'}\right)\right) \leq 9K_2. \tag{28}$$

*A.8.2. Variance of the first term*

Now, we can prove that the variance of each of the 10 separated sums converges to 0. The computation is very similar for each term, we only detail the most difficult one, and the same method apply directly to the other terms.

We consider:

$$\mathbb{V}\left[\frac{1}{p(n-1)^2}\sum_k\sum_{i,j}\left(\bar{y}_i\bar{y}_j - \frac{1}{n}\lambda_{ij}\right)^2\right] = \frac{p^2}{(n-1)^4 n^6}\sum_{\substack{k_1,k_1',k_2,k_2'\\k_3,k_3',k_4,k_4'}}\frac{1}{p^4}\sum_{i,j,\ell,m}\mathrm{Cov}\left(\left(y_{ik_1}y_{jk_1'} - \lambda_{ij}\mathbf{1}_{k_1=k_1'}\right)\left(y_{ik_2}y_{jk_2'} - \lambda_{ij}\mathbf{1}_{k_2=k_2'}\right),\right.$$
$$\left.\left(y_{\ell k_3}y_{mk_3'} - \lambda_{\ell m}\mathbf{1}_{k_3=k_3'}\right)\left(y_{\ell k_4}y_{mk_4'} - \lambda_{\ell m}\mathbf{1}_{k_4=k_4'}\right)\right).$$
(29)

If $k_1, k_1', k_2, k_2', k_3, k_3', k_4, k_4'$ are all different, then, the covariance is null. So there are at most $n^8 - n(n-1)(n-2)(n-3)(n-4)(n-5)(n-6)(n-7)$ non-zero terms in the sum over $k$-indices.

And we have from the preliminary inequality,

$$\sum_{i,j,\ell,m}\frac{1}{p^4}\mathrm{Cov}\left(\left(y_{ik_1}y_{jk_1'} - \lambda_{ij}\mathbf{1}_{k_1=k_1'}\right)\left(y_{ik_2}y_{jk_2'} - \lambda_{ij}\mathbf{1}_{k_2=k_2'}\right), \left(y_{\ell k_3}y_{mk_3'} - \lambda_{\ell m}\mathbf{1}_{k_3=k_3'}\right)\left(y_{\ell k_4}y_{mk_4'} - \lambda_{\ell m}\mathbf{1}_{k_4=k_4'}\right)\right) \leq 9K_2. \tag{30}$$

So,

$$\mathbb{V}\left[\frac{1}{p(n-1)^2}\sum_k\sum_{i,j}\left(\bar{y}_i\bar{y}_j - \frac{1}{n}\lambda_{ij}\right)^2\right] \leq \frac{p^2\left(n^8 - n(n-1)(n-2)(n-3)(n-4)(n-5)(n-6)(n-7)\right)}{(n-1)^4 n^6} \times 9K_2 \xrightarrow[n\to\infty]{} 0. \tag{31}$$

With the same method, we prove that the variance of the 9 other terms converges to 0 too. So, the first term has its variance converging to 0:

$$\mathbb{V}\left[\frac{1}{n^2}\sum_k\left\|\frac{n}{n-1}\tilde{x}_{\cdot k}\tilde{x}_{\cdot k}^\top - \Sigma\right\|^2\right] \xrightarrow[n\to\infty]{} 0. \tag{32}$$

*A.8.3. Variance of the second term*

The variance of the second term converges to 0 as $n$ goes to infinity. Indeed, from the proof of Lemma 3.4 in (Ledoit and Wolf, 2004) [10]:

$$\left(\frac{1}{n^2}\sum_k\left\|\frac{n}{n-1}\tilde{x}_{\cdot k}\tilde{x}_{\cdot k}^\top - S\right\|^2\right) - \left(\frac{1}{n^2}\sum_k\left\|\frac{n}{n-1}\tilde{x}_{\cdot k}\tilde{x}_{\cdot k}^\top - \Sigma\right\|^2\right) = \frac{1}{n}\|S - \Sigma\|^2. \tag{33}$$

And we have:

$$\mathbb{E}\left[\|S - \Sigma\|^4\right] \leq \mathbb{E}\left[(\|S - \mu I\| + \|\mu I - \Sigma\|)^4\right]. \tag{34}$$

$\mathbb{E}[\|S - \mu I\|^4]$ and $\|\mu I - \Sigma\|$ are bounded, from Lemmas 2 and 1 respectively, so $\mathbb{E}\left[\|S - \Sigma\|^4\right]$ is bounded. Consequently, $\mathbb{V}\left[\frac{1}{n}\|S - \Sigma\|^2\right] \underset{n\to\infty}{\longrightarrow} 0$. To conclude, $\mathbb{V}\left[\bar{b}^2\right] \underset{n\to\infty}{\longrightarrow} 0$.

### A.9. Proof of Lemma 7

Fully developing the variance and computing the expectations, we have:

$$\mathbb{V}[m] = \frac{1}{p^2(n-1)^2} \sum_{i,j} \left[\frac{n^2 - 2n + 1}{n} \mathbb{V}[y_{i1}y_{j1}] + \frac{n^2 - 1}{n} \lambda_{ij}^2 - \frac{n^2 - 2n + 1}{n} \lambda_{ii}\lambda_{jj}\right]$$

$$\mathbb{V}[m] = \frac{1}{p^2 n} \sum_{i,j} \mathbb{V}[y_{i1}y_{j1}] + \frac{n+1}{pn(n-1)}(\delta^2 - \beta^2 + \mu^2) - \frac{1}{n}\mu^2. \tag{35}$$

From the proof Lemma 1, we have:

$$\beta^2 = \frac{1}{pn} \sum_{i,j} \mathbb{V}[y_{i1}y_{j1}] + \frac{p+1}{n(n-1)}\mu^2 + \frac{1}{n(n-1)}(\delta^2 - \beta^2). \tag{36}$$

So, combining those two equations, we obtain:

$$\mathbb{V}[m] = \frac{1}{p}\left(\beta^2 - \frac{p+1}{n(n-1)}\mu^2 - \frac{1}{n(n-1)}(\delta^2 - \beta^2)\right) + \frac{n+1}{pn(n-1)}(\delta^2 - \beta^2 + \mu^2) - \frac{1}{n}\mu^2$$

$$\mathbb{V}[m] = q_0\beta^2 + q_1\delta^2 - q_2\mu^2, \tag{37}$$

which concludes the proof.

### A.10. Proof of Lemma 8

From the proof of Lemma 5, we have: $\mathbb{E}[m^2] = \mu^2 + \mathbb{V}[m]$, $\mathbb{E}[d^2] = \delta^2 - \mathbb{V}[m]$. And from Lemma 7, we have $\mathbb{V}[m] = q_0\beta^2 + q_1\delta^2 - q_2\mu^2$, which, when combined, immediately finishes the proof.

### A.11. Proof of Lemma 9

Combining Lemmas 5 and 8, we have:

$$\mathbb{E}[\bar{b}^2] = c_0\beta^2 + c_1\mathbb{E}[d^2] + c_2\mathbb{E}[m^2] + \frac{c_1 - c_2}{1 - q_1 - q_2}(q_0\beta^2 + q_1\mathbb{E}[d^2] - q_2\mathbb{E}[m^2])$$

$$\mathbb{E}[\bar{b}^2] = c_0^f \beta^2 + c_1^f \mathbb{E}[d^2] + c_2^f \mathbb{E}[m^2]. \tag{38}$$

Then, we deduce:

$$\mathbb{E}[b^2] = \frac{1}{c_0^f}\left((c_0^f\beta^2 + c_1^f\mathbb{E}[d^2] + c_2^f\mathbb{E}[m^2]) - c_1^f\mathbb{E}[d^2] - c_2^f\mathbb{E}[m^2]\right) = \beta^2. \tag{39}$$

So, $b^2$ is an unbiased estimator of $\beta^2$.

Concerning the quadratic mean convergence, we use the fact that the variances of $m^2$, $d^2$ and $\bar{b}^2$ converge to 0 as $n$ goes to infinity, from Corollary 2, and Lemmas 3 and 6 respectively, which concludes the proof.

### A.12. Proof of Lemma 10

For the upper and lower bounds:

$$b_u^2 - \beta^2 = \min(b_+^2, d^2) - \beta^2 \leq b_+^2 - \beta^2 \leq |b^2 - \beta^2| \leq \max(|b^2 - \beta^2|, |d^2 - \delta^2|)$$

$$b_u^2 - \beta^2 = \min(b_+^2 - \beta^2, d^2 - \beta^2) \geq \min(b_+^2 - \beta^2, d^2 - \delta^2) \geq -\max(|b^2 - \beta^2|, |d^2 - \delta^2|). \tag{40}$$

So,

$$\mathbb{E}[(b_u^2 - \beta^2)^2] \leq \mathbb{E}[\max(|b^2 - \beta^2|, |d^2 - \delta^2|)^2] \leq \mathbb{E}[(b^2 - \beta^2)^2] + \mathbb{E}[(d^2 - \delta^2)^2]. \tag{41}$$

Lemmas 3 and 9 trivially conclude the proof.

### A.13. Proof of Theorem 2

We will use the proof of Theorem 3.2 in (Ledoit and Wolf, 2004) [10] to prove ours. We check that we have the set of hypotheses required by the proof to work:

- $\alpha^2, \beta^2, \delta^2$ are non-negative, bounded, and $\alpha^2 + \beta^2 = \delta^2$
- $m - \mu$ converges to 0 in quartic mean,
- $d_u^2$ is nonnegative, and $d_u^2 - \delta^2 \xrightarrow[q.m]{} 0$,
- $0 \leq a_u^2 \leq d_u^2$ and $a_u^2 - \alpha^2 \xrightarrow[q.m]{} 0$,
- $a_u^2 + b_u^2 = d_u^2$, with $b_u^2 \geq 0$.

Then, we can apply the result of the theorem 3.2 from (Ledoit and Wolf, 2004), so $\mathbb{E}[\|S_n^* - \Sigma_n^*\|^2] \to 0$ and $S_n^*$ has the same asymptotic expected loss as $\Sigma_n^*$, i.e. $\mathbb{E}[\|S_n^* - \Sigma_n\|_n^2] - \mathbb{E}[\|\Sigma_n^* - \Sigma_n\|_n^2] \to 0$.

### A.14. Proof of Lemma 11

The proof from Ledoit and Wolf [10] can be applied here, as the hypotheses of their Lemma A.1 [10] are verified for $u^2 = |a_u^2 b_u^2 \delta^2 - \alpha^2 \beta^2 d_u^2|$, $\tau_1 = 2$ and $\tau_2 = 2$, from the same arguments they used.

### A.15. Proof of Theorem 3

Solving the convex minimization problem, we prove easily that $\Sigma^{**} = \mu I + \frac{\alpha_2}{d^2}(S - mI)$, with $\alpha_2 = \langle S, \Sigma \rangle - m\mu$, is a minimizer. We have then:

$$\|S^* - \Sigma^{**}\|^2 = (m - \mu)^2 + \frac{(a_u^2 - \alpha_2)^2}{d^2}. \tag{42}$$

$(m - \mu)^2$ converges to 0 in quadratic mean by Lemma 2. For the second term, we will use the Lemma A.1 from Ledoit and Wolf [10] with $u^2 = (a_u^2 - \alpha_2)^2$, $\tau_1 = 2$ and $\tau_2 = 0$. In the following, we check the assumptions of the Lemma A.1, i.e. $\frac{(a_u^2 - \alpha_2)^2}{d^2} \leq 2d^2 + 2\delta^2$ and $\mathbb{E}[(a_u^2 - \alpha_2)^2] \to 0$. We notice: $|\alpha_2| = |\langle \Sigma - \mu I, S - mI \rangle| \leq \|\Sigma - \mu I\| \|S - mI\| = \delta d$. It comes that:

$$\frac{(a_u^2 - \alpha_2)^2}{d^2} = \frac{a_u^4 + \alpha_2^2 - 2a_u^2 \alpha_2}{d^2} \leq \frac{2a_u^4 + 2\alpha_2^2}{d^2} \leq 2d^2 + 2\delta^2. \tag{43}$$

In order to prove $\mathbb{E}[(a_u^2 - \alpha_2)^2] \to 0$, let us show that $\mathbb{E}[(\alpha_2 - \alpha^2)^2] = \mathbb{V}[\alpha_2] \to 0$.

$$\mathbb{V}[\alpha_2] = \mathbb{V}[\langle S, \Sigma \rangle - m\mu] \leq 2\mathbb{V}[\langle S, \Sigma \rangle] + 2\mu^2 \mathbb{V}[m]. \tag{44}$$

$\mu$ is bounded by Lemma 1 and $\mathbb{V}[m] \to 0$ by Lemma 2, so $\mu^2 \mathbb{V}[m] \to 0$. Considering the other term, we have when developing:

$$\mathbb{V}[\langle S, \Sigma \rangle] = \frac{1}{n} \mathbb{V}\left[\frac{1}{p}\sum_i \lambda_{ii} y_{i1}^2\right] + \frac{2}{n(n-1)} \mathbb{V}\left[\frac{1}{p}\sum_i \lambda_{ii} y_{i1} y_{i2}\right] \leq \frac{n+1}{n(n-1)} K_2. \tag{45}$$

So, $\mathbb{V}[\langle S, \Sigma \rangle] \to 0$, and so $\mathbb{E}[(\alpha_2 - \alpha^2)^2] \to 0$. As $\mathbb{E}[(a_u^2 - \alpha^2)^2] \to 0$ by Lemma 9, it comes that $\mathbb{E}[(a_u^2 - \alpha_2)^2] \to 0$. Therefore, the assumptions of Lemma A.1 of (Ledoit and Wolf, 2004) [10] are verified by $u^2 = (a_u^2 - \alpha_2)^2$, $\tau_1 = 2$ and $\tau_2 = 0$. It proves that:

$$\mathbb{E}\left[\frac{(a_u^2 - \alpha_2)^2}{d^2}\right] \to 0. \tag{46}$$

Backing up, we have shown that $\mathbb{E}[\|S^* - \Sigma^{**}\|^2] \to 0$. We complete the proof of the theorem with the following inequality:

$$\mathbb{E}\left[\left|\|S^* - \Sigma\|^2 - \|\Sigma^{**} - \Sigma\|^2\right|\right] = \mathbb{E}\left[\left|\langle S^* - \Sigma^{**}, S^* + \Sigma^{**} - 2\Sigma \rangle\right|\right] \leq \sqrt{\mathbb{E}\left[\|S^* - \Sigma^{**}\|^2\right]} \sqrt{\mathbb{E}\left[\|S^* + \Sigma^{**} - 2\Sigma\|^2\right]}. \tag{47}$$

The first term converges to 0 as we showed above, and the second term is bounded because $\mathbb{E}[\|S^* - \Sigma\|^2]$ is bounded. So, the product converges to 0, which completes the proof.

### A.16. Proof of Theorem 4

The proof from Ledoit and Wolf [10] can be applied as it is, because it uses only the results of Theorem 3 which are the same as Theorem 3.3 in [10].

*A.17. Proof of Theorem 5*

Respectively from Lemmas 2 and 3, $m_r - \mu$ converges to 0 in quartic mean and $d_r^2 - \delta^2$ converges to 0 in quadratic mean. Let us define:

$$\bar{b}_r^2 = \frac{1}{(n-1)^2} \sum_{k=1}^{n} \|\tilde{x}_{\cdot k}(\tilde{x}_{\cdot k})^t - S\|^2 = \bar{b}^2 + \frac{1}{n(n-1)^2}\|S\|^2.$$

$\bar{b}^2 - \beta^2$ converges to 0 in quadratic mean by Lemma 6. Moreover, $\|S\|^2 - \mathbb{E}[\|S\|^2]$ converges to 0 in quadratic mean by Lemma 3, and $\mathbb{E}[\|S\|^2] = \beta^2 + \alpha^2 + \mu^2$, so $\mathbb{E}[\|S\|^2]$ is bounded by Lemma 1. So, $\frac{1}{n(n-1)^2}\|S\|^2$ converges to 0 in quadratic mean, which finally proves that $\bar{b}_r^2 - \beta^2$ converges to 0 in quadratic mean.

Following the idea of proof of Lemma 10, we have:

$$-\max(|\bar{b}_r^2 - \beta^2|, |d_r^2 - \delta^2|) \leq b_r^2 - \beta^2 \leq \max(|\bar{b}_r^2 - \beta^2|, |d_r^2 - \delta^2|). \tag{48}$$

Which, as previously, leads to the 2 following results: $\mathbb{E}[(b_r^2 - \beta^2)^2] \to 0$, $\mathbb{E}[(a_r^2 - \alpha^2)^2] \to 0$. We respect the set of hypotheses required by the proof of Theorem 2 to work, so identically we have that $\mathbb{E}[\|S_r^* - \Sigma^*\|^2] \to 0$ and $S_r^*$ has the same asymptotic expected loss as $\Sigma^*$, i.e. $\mathbb{E}\left[\left|\|S_r^* - \Sigma\|_n^2 - \|\Sigma^{**} - \Sigma\|^2\right|\right] \to 0$.

*A.18. Proof of Theorem 6*

Respectively from Lemmas 2 and 3, $m_m - \mu$ converges to 0 in quartic mean and $d_m^2 - \delta^2$ converges to 0 in quadratic mean. $\bar{b}^2 - \beta^2$ converges to 0 in quadratic mean by Lemma 6, and following the idea of proof of Lemma 10, we have:

$$-\max(|\bar{b}^2 - \beta^2|, |d_m^2 - \delta^2|) \leq b_m^2 - \beta^2 \leq \max(|\bar{b}^2 - \beta^2|, |d_m^2 - \delta^2|). \tag{49}$$

Which, as previously, leads to the 2 following results: $\mathbb{E}[(b_m^2 - \beta^2)^2] \to 0$, $\mathbb{E}[(a_m^2 - \alpha^2)^2] \to 0$. We respect the set of hypotheses required by the proof of Theorem 2 to work, so identically we have that $\mathbb{E}[\|S_m^* - \Sigma^*\|^2] \to 0$ and $S_m^*$ has the same asymptotic expected loss as $\Sigma^*$, i.e. $\mathbb{E}\left[\left|\|S_m^* - \Sigma\|^2 - \|\Sigma^* - \Sigma\|^2\right|\right] \to 0$.

*A.19. Proof of Theorem 7*

As $d_s^2 = d^2$ and $b_s^2 = b_m^2$, we trivially have the quadratic convergence to $\delta^2$ and $\beta^2$ respectively with Lemma 3 and Theorem 6. $m_s = \frac{n-1}{n}m$, and $m - \mu$ converges in quartic mean to 0 with $\mu$ bounded, so $m_s - \mu$ converges in quartic mean to 0. Similarly, $a_s^2 = \frac{n-1}{n}a_m^2$ and $a_m^2 - \alpha^2$ converges in quadratic mean with $\alpha^2$ bounded, so $a_s^2 - \alpha^2$ converges in quadratic mean. Finally, as $S_s^* = \frac{n-1}{n}S_m^*$, we similarly have $\mathbb{E}[\|S_s^* - \Sigma^*\|^2] \to 0$ and $S_s^*$ has the same asymptotic expected loss as $\Sigma^*$, i.e. $\mathbb{E}\left[\left|\|S_s^* - \Sigma\|^2 - \|\Sigma^* - \Sigma\|^2\right|\right] \to 0$.

*A.20. Proof of Lemma 12*

Let $X_{\cdot,k} \sim \mathcal{N}(0, \Sigma)$, $k \in [\![1, n]\!]$, $n$ iid samples. Trivially, $\mu = \frac{1}{p}\mathrm{tr}(\Sigma), \alpha^2 = \|\Sigma - \mu I\|^2$. For $\beta^2$, we use its development in the proof of Lemma 1.

$$\beta^2 = \frac{1}{pn}\sum_{i,j}\mathbb{V}[y_{i1}y_{j1}] + \frac{p+1}{n(n-1)}\mu^2 + \frac{1}{n(n-1)}\alpha^2. \tag{50}$$

As $X$ is Gaussian, we have for all $(i,j) \in [\![1, p]\!]^2$, $\mathbb{V}[y_{i1}y_{j1}] = \lambda_{ii}\lambda_{jj} + \lambda_{ij}^2$. So,

$$\frac{1}{pn}\sum_{i,j}\mathbb{V}[y_{i1}y_{j1}] = \frac{1}{n}\left(\|\Sigma\|^2 + \frac{1}{p}\mathrm{tr}(\Sigma)^2\right) = \frac{1}{n}(\alpha^2 + (p+1)\mu^2). \tag{51}$$

Which finally leads to,

$$\beta^2 = \frac{p+1}{n(n-1)}\mu^2 + \frac{1}{n(n-1)}\alpha^2 + \frac{1}{n}(\alpha^2 + (p+1)\mu^2) = \frac{p+1}{n-1}\mu^2 + \frac{1}{n-1}\alpha^2. \tag{52}$$

And, of course, $\delta^2 = \alpha^2 + \beta^2$.

*A.21. Proof of Lemma 13*

Let $\nu > 8$, $\Sigma$ a covariance matrix, $X_{\cdot,k} \sim t_\nu(0, \tilde{\Sigma})$, $k \in [\![1, n]\!]$, $n$ iid samples, with scale matrix $\tilde{\Sigma} = \frac{\nu-2}{\nu}\Sigma$. With this setup, we have as expected: $\mathbb{V}[X] = \Sigma$. Obviously, $\mu = \frac{1}{p}\mathrm{tr}(\Sigma), \alpha^2 = \|\Sigma - \mu I\|^2$.

For $\beta^2$, we will use the equation from the proof of Lemma 1, as in the Gaussian case.

$$\beta^2 = \frac{1}{pn}\sum_{i,j}\mathbb{V}[y_{i1}y_{j1}] + \frac{p+1}{n(n-1)}\mu^2 + \frac{1}{n(n-1)}\alpha^2. \tag{53}$$

To easily compute the variance term, we will use a characterization of multivariate $t$-distributions.

In fact, as for all $k \in [\![1, n]\!]$, $X_{\cdot,k} \sim t_\nu(0, \tilde{\Sigma})$, there exists two independent random variables $U_k$ and $Z_{\cdot,k}$ such that:

$$U_k \sim \chi_\nu^2, Z_{\cdot,k} \sim \mathcal{N}\left(0, \frac{\nu - 2}{\nu} \Lambda\right), Y_{\cdot,k} = \Sigma^{-1/2} X_{\cdot,k} = \sqrt{\frac{\nu}{U_k}} Z_{\cdot,k}. \tag{54}$$

Noticing that $\mathbb{E}\left[\frac{1}{U^2}\right] = \frac{1}{(\nu - 2)(\nu - 4)}$, we have,

$$\frac{1}{pn} \sum_{i,j} \mathbb{V}[y_{i1} y_{j1}] = \frac{1}{n}\left(\frac{\nu}{\nu - 4}(\alpha^2 + \mu^2) + \frac{\nu - 2}{\nu - 4} p \mu^2\right). \tag{55}$$

We can conclude,

$$\beta^2 = \left(\frac{\nu - 2}{(\nu - 4)n} + \frac{1}{n(n-1)}\right) p \mu^2 + \frac{1}{n}\left(\frac{\nu}{\nu - 4} + \frac{1}{n - 1}\right)(\alpha^2 + \mu^2). \tag{56}$$

And, of course, $\delta^2 = \alpha^2 + \beta^2$.